\documentclass[10pt]{article}
\usepackage{amsmath, amsfonts, amsthm, amssymb,verbatim}
  
        \newtheorem{theorem}{Theorem}[section]
\newtheorem{definition}[theorem]{Definition}
        \newtheorem{proposition}[theorem]{Proposition}
        \newtheorem{lemma}[theorem]{Lemma}
        \newtheorem{example}[theorem]{Example}
        \newtheorem{corollary}[theorem]{Corollary}
        \newtheorem{remark}[theorem]{Remark}

\numberwithin{equation}{section}

\newcommand \Bzero {{\mathbf B}_{\delta_0}}
\newcommand \Bone {{\mathbf B}_{\delta_1}}

\newcommand \Omegat {\widetilde \Omega}

\newcommand \bart    {\underline t}  
\newcommand \tbar    {\overline t}  
\newcommand \barT    {\underline T}  
\newcommand \Tbar    {\overline T}  
\newcommand \be     {\begin{equation}}
\newcommand \ee     {\end{equation}}
\newcommand \R    {\mathbb{R}} 
\newcommand \eps   {\epsilon} 
\newcommand \epsb  {\epsilon} 
\newcommand \Ccal   {\mathcal C}
\newcommand \Fcal   {\mathcal F}
\newcommand \Scal   {\mathcal S}
\newcommand \Lcal   {\mathcal L}
\newcommand \Rcal   {\mathcal R}
\newcommand \Hcal   {\mathcal H}
\newcommand \Ical   {\mathcal I}
\newcommand \Jcal   {\mathcal J}
\newcommand \del  {{\partial}}  
\newcommand \sgn  {\mbox{sgn}} 

\newcommand \ab    {\overline a} 
 
\newcommand \Ab    {{\overline {\!A}}}
  
\newcommand \wbf    {{\bold w}}

\newcommand \abf    {{\bold a}}  
\newcommand \Abf    {{\bold A}}  
  
\newcommand \Kcal    {\mathcal K}
\newcommand \Dcal    {\mathcal D} 
 
\newcommand \lam   {\lambda} 

\newcommand \lamb  {\overline\lambda}  
\newcommand \lb    {\overline l}  
\newcommand \rb    {\overline r} 
\newcommand \RN    {\R^N} 
  
\newcommand \wmin {w^{\min{}}}
\newcommand \wmax {w^{\max{}}}

\newcommand \cb {\overline c} 
\newcommand \eb {\overline e} 
 \newcommand \Osc {\text{osc}}

\begin{document} 

\title{Haar method, averaged matrix, wave cancellations, and $L^1$ stability for hyperbolic systems}

\author{Philippe G. LeFloch$^1$}

\footnotetext[1]{Laboratoire Jacques-Louis Lions \& Centre National de la Recherche Scientifique (CNRS),
Universit\'e Pierre et Marie Curie (Paris 6), 4 Place Jussieu,  75252 Paris, France. 
\\ 
E-mail : {pgLeFloch@gmail.com} 
                                                        \\ 
{\tt To cite this article:} J. Hyperbolic Differ. Equ. 3 (2006), 701--739. \hfill 
}

\renewcommand  \today {} 
\date{\today}
\maketitle

\begin{abstract}
We develop a version of Haar and Holmgren methods which applies to discontinuous solutions of 
nonlinear hyperbolic systems and allows us to control the $L^1$ distance between two entropy solutions. 
The main difficulty is to cope with linear hyperbolic systems with discontinuous coefficients.
Our main observation is that, while entropy solutions contain compressive shocks only,
the averaged matrix associated with two such solutions has compressive or undercompressive shocks, 
but no rarefaction-shocks --~which are recognized as a source for non-uniqueness and instability. 
Our Haar-Holmgren-type method rests on the geometry associated with the averaged matrix
and takes into account  adjoint problems and wave cancellations along generalized characteristics. 
It extends the method proposed earlier by LeFloch et al. for genuinely nonlinear systems.  
In the present paper, we cover solutions with small total variation and a class of systems with general flux
that need not be genuinely nonlinear and includes for instance fluid dynamics equations.  
We prove that solutions generated by Glimm or front tracking schemes 
depend continuously in the $L^1$ norm upon their initial data, by exhibiting 
an $L^1$ functional controling the distance between two solutions. 
\end{abstract}

\newpage 


\section{Introduction}
\label{IN-0}

Our purpose is to develop a version of Haar and Holmgren methods which applies to discontinuous solutions of 
nonlinear hyperbolic systems and allows us to control the $L^1$ distance between two entropy solutions. 
The main difficulty is to cope with linear hyperbolic systems with discontinuous coefficients.
Our main observation is that, while entropy solutions contain compressive shocks only,
the averaged matrix $\Ab(u,u')$ (see below) associated with two such solutions
$u,u'$ has compressive or undercompressive shocks, 
but no rarefaction-shocks --~which are recognized as a source for non-uniqueness and instability. 
 hence, the absence of rarefaction-shocks provides us with 
the key argument in order to establish the $L^1$ stability of solutions.
Our Haar-Holmgren-type method rests on the geometry associated with the averaged matrix
and takes into account  adjoint problems and wave cancellations along generalized characteristics.  
It extends the method proposed earlier by LeFloch and collaborators for the class of genuinely nonlinear systems
 \cite{LeFloch-IMA1,LeFlochXin,HuLeFloch,GoatinLeFloch-scalar,GoatinLeFloch-system,LeFloch-book}. 

In this paper, we present the framework in the context of the Haar method and 
cover solutions with small total variation and a class of systems with general flux
that need not be genuinely nonlinear and includes for instance fluid dynamics equations.  
Our main contribution is a proof that solutions generated by Glimm or front tracking schemes 
depend continuously in the $L^1$ norm upon their initial data. This is achieved by exhibiting 
an $L^1$ functional controling the distance between two solutions. 
 
A completely different strategy is developed by Bressan and collaborators
(see for instance \cite{AnconaMarson2,BCP,BressanMarson} and the references therein) 
who rely on the concepts of generalized tangent vectors and shift differentials. 
This method leads to technical difficulties and requires 
either special properties of the hyperbolic systems,
strong regularity assumptions on the solutions, 
or sophisticated restarting procedure to keep the piecewise regularity of the (approximate) solutions.  

We discuss here the well-posedness theory for strictly hyperbolic systems of conservation laws 
\be
\del_t u + \del_x f(u) = 0, \qquad u=u(t,x) \in \RN,  \, (t,x) \in \R_+ \times \R, 
\label{1.1}   
\ee
with general flux $f: \RN \to \RN$. In the present paper, we are motivated by a theorem established via 
independent methods in Bressan, Liu, and Yang \cite{BressanLiuYang,LiuYang-systems}
and Hu and LeFloch \cite{HuLeFloch} for the class of genuinely nonlinear systems. Our main 
objective is to cover a class of hyperbolic systems \eqref{1.1} that need not be genuinely nonlinear
and to establish that any solutions $u,u'$ generated by the Glimm scheme or 
by front tracking satisfy the {\sl $L^1$ continuous dependence property}  
\be 
\| u'(t) - u(t) \|_{L^1(\R)}  
\lesssim \|u'(0) - u(0)\|_{L^1(\R)}, \qquad t \in \R_+.
\label{1.2} 
\ee

We consider solutions with small amplitude whose range is 
included in the open ball $\Bzero \subset \RN$ centered at the origin and with radius $\delta_0$.  
In \eqref{1.1}, the flux ${f: \Bzero \to \RN}$ is a smooth mapping and, for all $u$, the Jacobian matrix 
$Df(u)$ admits $N$ distinct eigenvalues $\lam_1(u) < \ldots < \lam_N(u)$ and, therefore, basis of 
left- and right-eigenvectors denoted by $l_j(u), r_j(u)$ ($1 \leq j \leq N$). The latter are normalized so that 
right-eigenvectors have unit norm, and the basis of left- and right-eigenvectors are dual to each other. 
We refer to this as the ``standard normalization".

It is well-known that solutions to \eqref{1.1} are generally discontinuous 
and must be sought in the sense of distributions, typically within the class of functions with bounded variation. 
For the sake of uniqueness, the Cauchy problem associated with the system of conservation laws \eqref{1.1} 
must be supplemented with an entropy condition, originally 
proposed by Lax \cite{Lax} for genuinely nonlinear systems (satisfying, by definition, ${\nabla \lambda_i \cdot r_i \neq 0}$) 
and later extended to general equations and systems by Oleinik, Wendroff, and Liu \cite{Oleinik, Wendroff,Liu1}. 
See the textbooks \cite{Lax,Smoller,Dafermos-book,Bressan2,LeFloch-book} for details.  

Pioneering work by Glimm \cite{Glimm} via the so-called random choice scheme
established the existence of solutions to the Cauchy problem 
for genuinely nonlinear systems \eqref{1.1} and bounded variation (BV) initial data ${u_0 : \R \to \Bzero}$: 
\be
\label{1.initial}
u(0,\cdot ) = u_0. 
\ee
In recent years, extensive research on nonconvex flux functions 
initiated \cite{Liu1,Liu2,HayesLeFloch,AnconaMarson1,AnconaMarson2}
led to an existence theory for general systems \cite{Bianchini,IguchiLeFloch,LiuYang-exist}. 
In addition, Bianchini and Bressan \cite{BianchiniBressan} established the existence and continuous dependence of solutions
to general systems \eqref{1.1} via the vanishing viscosity method; their technique does not cover solutions 
generated by Glimm's scheme or front tracking algorithms. Further important material on the existence theory can also 
be found in earlier works \cite{Risebro,Bressan1,Bressan2,BressanLeFloch1,BressanLeFloch2,HoldenRisebro,LeFloch-book} 
on genuinely nonlinear systems.

To establish the continuous dependence property \eqref{1.2} we follow \cite{HuLeFloch,LeFloch-book}
and introduce a matrix-valued map $\Ab=\Ab(u,u')$ such that   
\be
\label{matrix0}
\aligned 
& \Ab(u, u') \, (u' - u) = f(u') - f(u),  
\\
& \Ab(u,u') = \Ab(u',u), \quad u,u' \in \Bzero. 
\endaligned
\ee
We refer to $\Ab(u,u')$ as the {\sl averaged matrix} of $u,u'$ and, for instance, we could choose
\be
\label{matrix1}
\Ab(u,u') := \displaystyle \int_0^1 Df \bigl( (1-\theta) \, u + \theta \, u' \bigr) \, d\theta, 
\ee
although a different choice may be preferable for certain systems \eqref{1.1}, for instance to handle solutions with large amplitude. 
Since $f$ is strictly hyperbolic we can assume (for a smaller value of $\delta_0$, if necessary)    
that $\Ab(u,u')$ admits $N$ real and distinct eigenvalues denoted by 
$$
\lamb_1(u,u') < \ldots < \lamb_N(u,u'),
$$
together with left- and right-eigenvectors $\lb_j(u,u'), \rb_j(u,u')$.

If $u=u(t,x)$ and $u'=u'(t,x)$ are solutions to the nonlinear system \eqref{1.1}, 
then clearly $\psi : = u - u'$ is a solution to 
the {\sl linear} hyperbolic system  
\be
\label{1.5}
\del_t \psi + \del_x \bigl( \Abf \, \psi \bigr) = 0, 
\ee
where $\Abf := \Ab(u, u')$. Consequently, in order to establish \eqref{1.2} it is sufficient to derive the $L^1$ stability property 
\be
\| \psi(t) \|_{L^1(\R)}  
\lesssim  \|\psi(0)\|_{L^1(\R)}, \quad t \in \R_+, 
\label{1.6}
\ee 
for a sufficiently large class (covering the situation of interest) 
of matrix-valued mappings $\Abf=\Abf(t,x)$ and solutions $\psi=\psi(t,x)$ (which need not be of the form
$\Abf= \Ab(u, u')$ and $\psi = u-u'$).  

One should realize that the stability estimate \eqref{1.6} does not follow from the standard theory 
of linear hyperbolic equations, since \eqref{1.5} has {\sl discontinuous} coefficients. In fact,
we {\sl should not} expect the estimate \eqref{1.6} to be valid for {\sl arbitrary} systems \eqref{1.5}; see Example~\ref{Example1}. 
One of our objectives is precisely to identify a sufficiently large class of systems \eqref{1.6} covering
the situation of interest and, especially, 
to identify properties that the averaged matrix $\Abf(u,u')$ must satisfy.

\
 
In most of the present paper, we assume that $\psi$ and $\Abf$ are of the form 
$\psi = u - u'$ and $\Abf := \Ab(u, u')$ and postpone the discussion of more general solutions and matrices. 
In practice, we will obtain \eqref{1.2} (or \eqref{1.6}) only in the limit of 
a sequence of approximate solutions $u_h$, $u_h'$: 
\be 
\| u_h(t) - u_h'(t) \|_{L^1(\R)}  
\lesssim \|u_h(0) - u_h'(0)\|_{L^1(\R)} + o(1), 
\quad t \geq 0, 
\label{1.30} 
\ee
where $o(1)$ tends to zero with the discretization parameter $h$.

The approach above was used by Haar \cite{Haar} to establish 
uniqueness results for {\sl smooth} solutions of partial differential equations. 
It is also formally equivalent (by duality) to Holmgren method. 
The method was applied to nonlinear hyperbolic equations with special structure 
by several authors, beginning with Oleinik \cite{Oleinik} for scalar equations.
(Cf.~\cite{Dafermos-book,LeFlochXin} and the references therein.)  
Using this strategy to handle {\sl discontinuous} solutions of {\sl general nonlinear} 
systems \eqref{1.1} was first proposed by LeFloch 
and successfully led to a proof of \eqref{1.2} for genuinely nonlinear systems 
\cite{HuLeFloch}. Independently, a method based on an explicit weighted $L^1$ distance 
was introduced by Liu and Yang \cite{LiuYang-scalar} for scalar conservation laws with convex flux 
and later extended to genuinely nonlinear systems by Bressan, Liu, and Yang \cite{BressanLiuYang,LiuYang-systems}.

\

The content of the present paper is as follows. 
\begin{enumerate} 
\item Recall that a rarefaction-shock in the $i$-family, by definition, 
has all $i$-characteristics moving away from it. We prove that rarefaction-shocks are the only source of 
non-uniqueness of solutions to \eqref{1.5} and provide an $L^1$ estimate for piecewise constant solutions of \eqref{1.5}: 
$$
\| \psi(t)\|_{L^1} \lesssim \| \psi(0)\|_{L^1} + \Dcal_{\Rcal^\Abf}(t),  \qquad t \in \R_+, 
$$ 
where the error term $\Dcal_{\Rcal^\Abf}(t)$ involves only rarefaction shocks in the matrix $\Abf$.  
\item Generalizing a similar observation made earlier for genuinely nonlinear systems \eqref{1.1},  
we prove that the averaged matrix $\Ab(u,u')$ associated with two entropy solutions $u,u'$ 
does not admit rarefaction-shocks. More precisely, this property holds 
for ``sufficiently robust'' shocks, in a sense defined below (Section~\ref{section3}). 
Together with our $L^1$ stability estimate, this gives a strong indication that \eqref{1.6} should hold. 
\item Next, we actually construct an $L^1$ functional satisfying the condition required for stability. 
Our new functional can be regarded as a generalization of functionals proposed independently by 
Liu and Yang (via an explicit formula) 
and LeFloch (via a constructive method) for genuinely nonlinear systems. 
The main novelty is that not only the sign of (the characteristic components of) $u' - u$ must be taken 
into account but also the {\sl monotonicity} of the waves (the sign of $u_+ - u_-$ or $u_+' - u_-'$ 
at jumps) and the {\sl jump of the averaged speed} $\lamb_i(u_+,u_+') - \lamb_i(u_-,u_-')$.
Indeed, the {\sl eigenvalues} of the averaged matrix $\Ab(u,u')$ and their jumps across shocks
are found to play a central role (instead of the usual wave strengths).  
\item Furthermore, our construction of the weight takes into account wave cancellations, 
by solving a backward evolution problem for an adjoint equation associated with \eqref{1.5}. 
Wave partitions (Liu \cite{Liu2}) are essential to exhibit the cancellation effect. 
Our present approach is in fact a generalization of the construction of the weight proposed in \cite{HuLeFloch,LeFloch-book}
for of genuinely nonlinear systems. (Namely, 
trains of oscillating waves with alternating signs and equal strength already arose in  \cite{HuLeFloch,LeFloch-book}.) 
\item We conclude with the $L^1$ continuous dependence property \eqref{1.2} for a large class of hyperbolic systems, 
including the equations arising in isentropic fluid dynamics. 
\end{enumerate} 

Our method should be useful in other applications when wave cancellations allow to derive sharper bounds on solutions. 
For the Lagrangian and Eulerian compressible fluid equations, our method applies to solutions with large total variation. 
We refer to a follow-up paper for further developments. 

 
\section{Shock waves and characteristic components}
\label{LI-0}
 
\subsection{Compressive, undercompressive, and rarefaction shocks} 

We will rely on the following terminology and notation. In this section and in Section~\ref{section4} 
the maps ${\psi=\psi(t,x)}$ and ${\Abf=\Abf(t,x)}$ are piecewise constant
and admit finitely many polygonal lines of discontinuity, which can intersect two at a time and 
at finitely many times only; we refer to these lines as shock waves. Consider, for instance, the vector-valued map $\Abf$,
together with the associated decomposition 
$$ 
\R_+ \times \R = \Ccal^\Abf \cup \Jcal^\Abf \cup \Ical^\Abf, 
$$
where $\Ccal^\Abf$ is the set of points where $\Abf$ is locally constant,
$\Jcal^\Abf$ the set of points where $\Abf$ has a shock, and $\Ical^\Abf$ the set of points where two shock waves meet. 
 Also, $\Jcal^\Abf(t)$ will denote the slice of $\Jcal^\Abf$ 
consisting of jump points with time component $t$. 
At each jump point $(t,x) \in \Jcal^\Abf$, denote the left- and right-hand limits of $\Abf$ by $\Abf^\pm = \Abf^\pm(t,x)$
and the shock speed by $\lamb^\Abf=\lamb^\Abf(t,x)$. 
 
We assume that, at each $(t,x) \in \Ccal^\Abf$, the matrix $\Abf(t,x)$ admits $N$ real and distinct eigenvalues 
$$
\lam_1^\Abf(t,x) < \ldots < \lam_N^\Abf(t,x)
$$  
and, therefore, left- and right-eigenvectors $l_j^\Abf(t,x), r_j^\Abf(t,x)$ ($1 \leq j \leq N$), normalized 
in the standard way (see introduction). 
We always assume that the basis $l_j^\Abf$ and $r_j^\Abf$ remain sufficiently close to fixed basis of constant vectors. 
 At each $(t,x) \in \Jcal^\Abf$, the notation $l_{j\pm}^\Abf = l_{j\pm}^\Abf(t,x)$
and $r_{j\pm}^\Abf = r_{j\pm}^\Abf(t,x)$ for the corresponding left- and right-hand limits will be used. 

We also assume that the eigenvalues are uniformly separated, in the sense
that there exist fixed and disjoint intervals  $\bigl[\lam_j^{\min}, \lam_j^{\max}\bigr]$ ($j=1, \ldots, N$) 
so that at each $(t,x) \in \Ccal^\Abf$
$$
\lam_j^{\min} \leq \lam_j^\Abf(t,x) \leq \lam_j^{\max},  
$$ 
while for each $(t,x) \in \Jcal^\Abf$ there exists an index $i$ such that  
$$
\lam_i^{\min} \leq \lamb^\Abf(t,x) \leq \lam_i^{\max}.  
$$  
This allows us to decompose the set of shocks in $\Abf$, in form 
$\Jcal^\Abf = \Jcal_1^\Abf \cup \ldots \cup \Jcal_N^\Abf$, where $\Jcal_i^\Abf$ is the set of all $i$-shocks.

Our results will be uniform for all matrices satisfying the uniform total variation bound
\be
\label{BV}
TV(A(t)) \leq C, \qquad t \geq 0, 
\ee
for some fixed constant $C>0$. 
A matrix-valued field $\Abf=\Abf(t,x)$ fulfilling all of the conditions above will be called {\sl uniformly hyperbolic
with bounded variation,} and we will be interested in estimates that are valid for the whole class of such matrix-fields.

\begin{definition} 
\label{def21}
An $i$-shock wave in $\Jcal_i^\Abf$ is called 
\begin{enumerate}
\item a {\rm compressive (or Lax)} shock if $\lam_{i-}^\Abf \geq \lamb^\Abf \geq \lam_{i+}^\Abf$, 
\item a {\rm slow undercompressive} shock if $\min\bigl(\lam_{i-}^\Abf, \lam_{i+}^\Abf\bigr) \geq \lamb^\Abf$, 
\item a {\rm fast undercompressive} shock if  $\max\bigl(\lam_{i-}^\Abf, \lam_{i+}^\Abf\bigr) \leq \lamb^\Abf$,
\item or a {\rm rarefaction} shock if $\lam_{i-}^\Abf < \lamb^\Abf < \lam_{i+}^\Abf$.  
\end{enumerate} 
\end{definition}

Based on this definition, a shock which, for instance, is characteristic on its left-hand side and 
compressive on its right-hand side 
can be regarded as either a compressive or a fast undercompressive shock. 
For convenience in certain formulas, we may also allow large inequalities in the definition of rarefaction-shocks. 
 
The notation $\Lcal^\Abf$, $\Scal^\Abf$, $\Fcal^\Abf$, and $\Rcal^\Abf$ will be used 
for the sets of all compressive, slow undercompressive, fast undercompressive, and rarefaction shocks, respectively, 
and $\Lcal_i^\Abf$, $\Scal_i^\Abf$, $\Fcal_i^\Abf$, and $\Rcal_i^\Abf$ will denote the subsets associated with 
the $i$-characteristic family. Hence,  
$$
\Jcal^\Abf_i = \Lcal_i^\Abf \cup \Scal_i^\Abf \cup \Fcal_i^\Abf \cup \Rcal_i^\Abf. 
$$

\begin{example} 
\label{Example1} 
It is not difficult to see that when $N=1$ and $\Abf = \abf : \R_+\times\R \to \R$ is chosen to be 
$$
\abf(t,x) = \begin{cases}
 -1,  &   x<0, 
\\
\hskip.23cm 1,  &   x>0, 
\end{cases}
$$
then the Cauchy problem associated with the scalar equation 
\be
\label{2.3b}
\del_t \psi + \del_x (\abf \, \psi ) = 0, \qquad \psi=\psi(t,x) \in \R  
\ee
admits infinitely many solutions.  Therefore, if the speed coefficient $\abf$ admits rarefaction-shocks then
the (uniqueness property and, a~fortiori, the) continuous dependence property \eqref{1.6} can not hold
for solutions of \eqref{2.3b}.  
\end{example}
 

To establish \eqref{1.6} one possible strategy 
is to determine a weighted $L^1$ norm that is decreasing in time 
for all solutions of \eqref{1.5}. In view of Example~\ref{Example1} above, we will need to restrict attention
to matrix $\Abf$ with no rarefaction-shocks. 
Introduce the {\sl characteristic components} $\alpha = (\alpha_1, \ldots, \alpha_N)$ of a solution $\psi$, defined by 
\be
\label{charact}
\psi(t,x) =: \sum_j \alpha_j(t,x) \, r_j^\Abf(t,x), \qquad (t,x) \in \Ccal^\Abf \cap \Ccal^\psi. 
\ee
All the summations are over $\{1, \ldots, N\}$, except if indicated otherwise.  Then, 
consider the weighted norms
$$ 
\| \psi(t)\|_{\wbf(t)}  := \int_\R \sum_j |\alpha_j(t,x)| \, \wbf_j(t,x) \, dx,  
$$ 
where the piecewise constant weight $(\wbf_1, \ldots, \wbf_N)$ satisfies
\be
\label{LI-bounds}
\wmin \leq \wbf_j(t,x) \leq \wmax
\ee 
for all $j$ and $(t,x)$. Clearly, for fixed constants $\wmin, \wmax>0$, 
these norms are equivalent to the standard $L^1$ norm. 

Following \cite{HuLeFloch,LeFloch-book}
we require that $\wbf$ be a solution of the adjoint system associated with \eqref{1.5} 
\be
\label{2.adjoint}
\del_t \wbf + \Abf \, \del_x \wbf = g,
\ee
where the source-term $g$ may consist of measure-terms and will be specified later
when we will impose additional constraints on the jumps of the weight.  
(Strictly speaking, at a shock of $\Abf$ a definition of the nonconservative product $\Abf \, \del_x \wbf$ is necessary, 
and indeed a suitable definition 
can be given along the same lines as in \cite{LeFloch1,LeFloch-IMA2,DLM,LeFlochLiu}. 
In the present paper we bypass the difficulty and need not address this issue directly.) 

Precisely, we assume that every shock in $\wbf$ coincides with either a shock in $\Abf$ or a characteristic straightline, 
so that for all but finitely many (interaction) points $(t,x) \in \Jcal^\wbf$ we have either 
\begin{enumerate}
\item $(t,x) \in \Jcal^\Abf$, or 
\item there exists an index $i$ such that $\lamb^\wbf(t,x) = \lam_i^\Abf(t,x)$ 
and, for all $j \neq i$, $\wbf_j^+(t,x) = \wbf_j^-(t,x)$. 
\end{enumerate}
In addition, we could also allow the components $\wbf_j$ to exhibit a decreasing jump in time, that is, 
\be
\label{decrea}
\wbf_j(t+,x) \leq \wbf_j(t-,x), \quad (t,x) \in \R_+ \times \R. 
\ee
We will say in short that $\wbf$ {\sl formally solves the adjoint system} \eqref{2.adjoint} and, later, 
will also impose further constraints on the jumps $\wbf_j^+ - \wbf_j^-$. 

The following two lemmas were established in \cite{HuLeFloch, LeFloch-book}. 

\begin{lemma}[Time-derivative of the weighted norm] 
\label{weight} 
Suppose that the weight $\wbf$ formally solves the adjoint system. Then, 
for any piecewise constant $\psi$ the weighted norm satisfies 
\be
\label{2.8}
\aligned 
{d \over dt} \| \psi(t) \|_{\wbf(t)}  
\leq & \sum_{i,j}  \sum_{x \in \Jcal_i^\Abf(t)} \beta_j^-(t,x) \, \wbf_j^-(t,x) + \beta_j^+(t,x) \, \wbf_j^+(t,x) 
\endaligned 
\ee
for all but finitely many $t$, where
\be
\label{2.defbeta}
\aligned 
& \beta_j^-(t,x) := \bigl(\lamb^\Abf(t,x) - \lam_{j-}^\Abf(t,x) \bigr) 
  \, |\alpha_j^-(t,x)|,   
\\ 
& \beta_j^+(t,x) := \bigl(\lam_{j+}^\Abf(t,x) - \lamb^\Abf(t,x)\bigr) 
  \, |\alpha_j^+(t,x)|. 
\endaligned
\ee
are referred to as the {\rm characteristic flux.} 
\end{lemma}

\begin{lemma} [Signs of the characteristic flux] 
\label{lemma24} 
Across each $i$-shock the components $\beta_j^\pm$ satisfy for all $j \neq i$ 
\be
\label{sign-j}
\aligned 
& \pm \beta_j^\pm \leq 0, \quad j<i, 
\\
& \pm \beta_j^\pm \geq 0, \quad j>i, 
\endaligned 
\ee
and for $j= i$ 
\be
\label{sign-i}
\aligned 
&     \hskip.35cm  \beta_i^\pm \leq 0,  \quad \Lcal_i^\Abf, 
       \qquad \qquad 
       & \beta_i^\pm \geq 0,  \quad \Rcal_i^\Abf,   
\\
&\pm \beta_i^\pm \geq 0,  \quad \Scal_i^\Abf,   
\qquad \qquad 
      & \pm \beta_i^\pm \leq 0,  \quad  \Fcal_i^\Abf. 
\endaligned
\ee
\end{lemma} 

Clearly, the signs (and amplitude) of the terms ${\lamb^\Abf(t,x) - \lam_{j\pm}^\Abf(t,x)}$
are essential in the evolution of the weighted norm. For $j \neq i$ one of these two terms 
always has a favorable sign. 
For $j =i$ the signs depend upon the nature of the shock, and it is only in the 
case of a rarefaction shock that none of these two terms have a favorable sign. 


\subsection{Properties of characteristic flux}

Our estimates below depend upon the strength of $i$-shocks in the matrix ${\Abf=\Abf(t,x)}$, measures 
with the {\sl total strength} $| \Abf^+(t,x) - \Abf^-(t,x) |$, or in term of the $i$-eigenvector $| r_{i+}^\Abf(t,x) - r_{i-}^\Abf(t,x) |$. 
It is convenient to introduce also the characteristic flux with the modulus suppressed in \eqref{2.defbeta}. 
  
\begin{lemma}[Jump relation for characteristic flux across an $i$-shock]
\label{jump}
Let $\psi=\psi(t,x)$ be a solution of the system \eqref{1.5}. 
Then, at every shock point $(t,x)$ the characteristic flux
$$
\gamma_j^\pm := \bigl(\lamb^\Abf - \lam_{j\pm}^\Abf\bigr) \, \alpha_j^\pm
$$
satisfy 
\be
\label{2.13}
\gamma_j^+ = \gamma_j^- + \sum_k  l_{j+}^\Abf \, \bigl(r_{k-}^\Abf - r_{k+}^\Abf\bigr) \, \gamma_k^-
\ee
and
\be
\label{2.13b}
 \gamma_j^+  =  \gamma_j^- + {r_{j-}^\Abf + r_{j+}^\Abf \over 1 + r_{j+}^\Abf  \cdot r_{j-}^\Abf} \cdot 
 \sum_{k \neq j} \big(\gamma_k^- \, ( r_{k-}^\Abf - r_{k+}^\Abf ) + (\gamma_k^- - \gamma_k^+) \, r_{k+}^\Abf \big). 
\ee
\end{lemma} 

Lemma~\ref{jump} implies the expansion 
\be
\label{beta}
\gamma_j^+ = \gamma_j^- + O(| r_{i+}^\Abf - r_{i-}^\Abf| ) \, |\gamma_i^-| 
                       + O(|\Abf^+ - \Abf^-|) \, \sum_{k\neq i} |\alpha_k^-|, \quad 1 \leq j \leq N,
\ee 
and for $j \neq i$ and $j=i$, respectively, 
\be
\label{beta2} 
\aligned 
& \gamma_j^+ = \gamma_j^- + O(1) \sum_{k \neq j} |\gamma_k^+ - \gamma_k^-| + O(| r_{i+}^\Abf - r_{i-}^\Abf| ) \, |\gamma_i^-|  
                 + O(|\Abf^+ - \Abf^-|) \, \sum_{k \neq i, j} |\alpha_k^-|,  
\\
& \gamma_i^+ = \gamma_i^- + O(1) \sum_{k \neq i} |\gamma_k^+ - \gamma_k^-| + O(|\Abf^+ - \Abf^-|) \, \sum_{k \neq i} |\alpha_k^-|. 
\endaligned  
\ee
Observe that \eqref{beta2} contain no error term in the indices $j$ and $i$, respectively. 
For instance, the second identity shows that
$$
\gamma_i^+ = \gamma_i^- + O(1) \sum_{k \neq i \atop \pm} |\alpha_k^\pm|, 
$$
and, therefore, if $\alpha_j^- = \alpha_j^+ =0$ for all $j\neq i$ 
then $\gamma_i^+ = \gamma_i^-$.  

\begin{proof} The Rankine-Hugoniot relation associated with \eqref{1.5} reads  
$$
- \lamb^\Abf \, \bigl(\psi_+ - \psi_-\bigr) + \Abf_+ \psi_+ - \Abf_- \psi_- = 0, 
$$ 
thus 
\be
\label{2.identity}
\aligned 
\sum_k  \gamma_k^+ \, r_{k+}^\Abf 
& = \sum_k \bigl(\lamb^\Abf - \lam_{k+}^\Abf\bigr) \, \alpha_k^+ \, r_{k+}^\Abf 
\\
& = 
\sum_k \bigl(\lamb^\Abf - \lam_{k-}^\Abf \bigr) \, \alpha_k^- \, r_{k-}^\Abf
= 
\sum_k  \gamma_k^- \, r_{k-}^\Abf.
\endaligned 
\ee
Multiplying this identity by each eigenvector $l_{j+}^\Abf$ and using the normalization 
$l_{j+}^\Abf \, r_{k+}^\Abf = 0$ if $j \neq k$ and $l_{j+}^\Abf \, r_{j+}^\Abf = 1$, 
we arrive at  \eqref{2.13}. 

For illustration, suppose for instance that, for some fixed $j$, all of the components 
$\gamma_k^\pm$ for $k \neq j$ vanish. Then, the identity \eqref{2.identity} implies that  
$$ 
\gamma_j^+ \, r_{j+}^\Abf = \gamma_j^- \, r_{j-}^\Abf. 
$$ 
Hence, either $\gamma_j^- = \gamma_j^+ = 0$, or else 
$r_{j-}^\Abf = r_{j+}^\Abf$ (since both are unit vectors and are close to each other) and, in both cases,    
$$ 
\gamma_j^+ = \gamma_j^-. 
$$ 
We observe that there should be no error term
(for instance in \eqref{beta}) when the components $\gamma_k^\pm$ vanish for all $k$ 
but one single index. 

We can reformulate the estimate \eqref{2.13} by observing that for all $j$ 
$$
\aligned 
\gamma_j^+ \, r_{j+}^\Abf -  \gamma_j^- \, r_{j-}^\Abf
& = \sum_{k \neq j}  \gamma_k^- \, r_{k-}^\Abf -  \gamma_k^- \, r_{k-}^\Abf 
\\
& = \sum_{k \neq j} \gamma_k^- \, ( r_{k-}^\Abf - r_{k+}^\Abf ) + (\gamma_k^- - \gamma_k^+) \, r_{k+}^\Abf =: \Omega. 
\endaligned 
$$ 
Since the $r_j^\pm$'s are unit vectors
\be
\label{system2}
\aligned 
& 
\gamma_j^+ = \gamma_j^- \, r_{j+}^\Abf  \cdot r_{j-}^\Abf + r_{j+}^\Abf \cdot \Omega, 
\\
& 
\gamma_j^- = \gamma_j^+ \, r_{j+}^\Abf \cdot r_{j-}^\Abf - r_{j-}^\Abf \cdot \Omega,
\endaligned 
\ee 
this implies 
$$
\big( 1 - (r_{j+}^\Abf  \cdot r_{j-}^\Abf)^2 \big) \, \gamma_j^+ =  \big(  r_{j+}^\Abf    - (r_{j+}^\Abf  \cdot r_{j-}^\Abf) r_{j-}^\Abf \big) \cdot \Omega,
$$
and, since $r_{j+}^\Abf  \cdot r_{j-}^\Abf = 1 + o(1)$, 
$$
\big( 1 - r_{j+}^\Abf  \cdot r_{j-}^\Abf \big) \, \gamma_j^+ 
= {\big(  r_{j+}^\Abf    - (r_{j+}^\Abf  \cdot r_{j-}^\Abf) r_{j-}^\Abf \big)  \over 1 + r_{j+}^\Abf  \cdot r_{j-}^\Abf} \cdot \Omega.
$$
Combining this estimate with the second identity in \eqref{system2} we obtain \eqref{2.13b}.
\end{proof} 
 

We fix a (sufficiently small) constant ${\kappa_1>0}$ and introduce the following terminology. 

\begin{definition} 
\label{domin}
Consider a solution $\psi= \psi(t,x)$ of \eqref{1.5} together with its characteristic flux
$\beta_j$ defined in \eqref{2.defbeta}. 
Then, at any $i$-shock the $j$-characteristic flux $(1 \leq j \leq N$) is said to be {\rm dominant} if 
\be
\label{dom-1}
\kappa_1 \, |\beta_j^-| \geq | r_{i+}^\Abf - r_{i-}^\Abf| \, |\beta_i^-| + |\Abf^+ - \Abf^-| \, \sum_k |\beta_k^-|. 
\ee 
\end{definition}

The notion of dominant flux is quite natural in view of the expansion \eqref{beta}.
Indeed, it follows that for all dominant characteristic flux 
\be
\label{signes}
\sgn(\gamma_j^+) = \sgn(\gamma_j^-),   \quad \text{ $j$ dominant,}  
\ee
which leads us to the following conclusion. 

\begin{lemma} [Signs of characteristic components]
For all $j \neq i$ 
\be
\label{signes-j}
\sgn(\alpha_j^+) = \sgn(\alpha_j^-),     \quad \text{ $j$ dominant,}  
\ee
while for $j=i$ 
\be
\label{signes-i}
\sgn(\alpha_i^+) = 
\begin{cases} 
\hskip.25cm \sgn(\alpha_i^-), \quad & \text{ $\Lcal_i \cup \Rcal_i$ and $i$ dominant,}  
\\
-\sgn(\alpha_i^-), \quad & \text{ $\Scal_i \cup \Fcal_i$ and $i$ dominant.}  
\end{cases} 
\ee
\end{lemma} 

\begin{remark} 
\label{domin-remark} 
1. From the expansion \eqref{beta} and the inequality \eqref{dom-1}, 
an analogous inequality follows with $\beta_j^-$ replaced by $\beta_j^+$, i.e. 
$$ 
2 \kappa_1 \, |\beta_j^+| \geq | r_{i+}^\Abf - r_{i-}^\Abf| \, |\beta_i^-| + |\Abf^+ - \Abf^-| \, \sum_k |\beta_k^-|. 
$$  
Conversely, if the latter inequality holds, then $\beta_j^-$ satisfies \eqref{dom-1} with $\kappa_1$ replaced by $3 \kappa_1$.  

2. At this stage of the discussion, in view of Lemma~\ref{lemma24} it should be clear that 
for the right-hand side of \eqref{2.8} to be non-negative
the weight should jump up or down at $i$-shocks, according to 
\be
\label{updown}
\wbf_j^+ - \wbf_j^- 
\begin{cases} 
\geq 0, &  j < i,
\\
\leq 0, & j = i \text{ and } \Scal_i,   
\\
\geq 0, & j = i \text{ and } \Fcal_i, 
\\
\leq 0, &  j > i. 
\end{cases} 
\ee
\end{remark}


\section{The averaged matrix admits no rarefaction-shocks}
\label{section3}

\subsection{Scalar and fluid dynamics equations} 

Our aim in the present section 
is to derive properties satisfied by the averaged matrix $\Ab(u,u')$ (see \eqref{matrix0}
and, for instance, \eqref{matrix1})
associated with two solutions $u,u'$ of \eqref{1.1}. 
Recall that a shock is {\sl entropy admissible} if it satisfies the standard entropy condition \cite{Oleinik,Lax,Wendroff,Liu1}; 
see ~\eqref{3.8} and \eqref{3.111} below. 
Our main observation is that entropy admissible shocks in weak solutions $u,u'$ 
can not be rarefaction-shocks for the matrix-field $\Ab(u,u')$. 
Throughout this section, we consider weak solutions with bounded variation, which are known
to admit, in a measure-theoretic sense, shock curves and left- and right-hand traces.

Precisely, we consider maps $u$
in the space $BV(\R_+ \times \R,\R^N) \cap L^\infty(\R_+ \times \R,\R^N)$
of (bounded) functions with bounded variation. To such $u$
we can associate a decomposition of $\R_+\times \R$ (e.g.~Federer \cite{Federer})
\be
\label{1}
\R_+ \times \R = \Ccal^u \cup \Jcal^u \cup \Ical^u,
\ee
where $\Ccal^u$ is the set of points of approximate continuity of $u$, $\Jcal^u$ the set of
approximate jump discontinuity, and $\Hcal^1(\Ical^u) = 0$.
The set $\Jcal^u$ is $H^1$-rectifiable (for the one-dimensional Hausdorff measure) 
and so, except for a set of $H^1$-measure zero, can be covered by a countable family of $C^1$ graphs. 
In the $L^1$ sense, we can also define
the approximate limits $u_\pm(t,x)$ from each side of the 
shock at $(t,x)$, and $\lamb^u(t,x)$ denotes the shock speed. 

\begin{lemma}[Characterization of the nature of shocks] 
\label{Lemma3.2} Consider a scalar conservation law with general flux $f:\R \to \R$
and a shock wave connecting the states $u_-, u_+$ 
and propagating at the Rankine-Hugoniot speed 
\be 
\ab(u_-, u_+) := \int_0^1 \del_u f(s \, u_- + (1- s) \, u_+) \, ds. 
\label{3.1}
\ee
For any real $u'$ the nature of the shock wave in the averaged speed function 
\be
\abf(t,x) = \begin{cases} 
\abf_- := \ab(u_-,u'),  &  x <  \ab(u_-, u_+) \, t, 
\\
\abf_+ := \ab(u_+,u'),  &  x >  \ab(u_-, u_+) \, t, 
\end{cases}
\label{3.2}
\ee 
is uniquely determined by the signs of $\rho(u_-,u_+,u') :=  (u_- - u') \, (u_+ - u_-)$
and the jump $(\abf_+ - \abf_-) $ of the averaged speed, 
as follows: 
\begin{enumerate} 
\item If $u_\pm - u'$ have the {\rm same sign}, 
\be 
\rho(u_-,u_+,u') \, (\abf_+ - \abf_-) 
\begin{cases} 
\leq 0, & \Scal^\abf,     
\\ 
\geq 0, & \Fcal^\abf. 
\end{cases}  
\label{3.3}  
\ee
\item If $u_\pm - u'$ have {\rm opposite signs}, 
\be
\rho(u_-,u_+,u') \, (\abf_+ - \abf_-) 
\begin{cases} 
\leq 0, & \Rcal^\abf,     
\\ 
\geq 0, & \Lcal^\abf. 
\end{cases}  
\label{3.4}  
\ee
\end{enumerate}
\end{lemma} 
 
With some abuse of notation, in the first case of \eqref{3.4} we have allowed here the case of equalities
in the definition of a rarefaction-shock. 
  
\begin{proof} The inequalities can be checked geometrically from the graph of $f$, or else
algebraically from the jump relation, as we now do. In view of 
$$
\aligned 
& \ab(u_-,u_+) \, (u_+ - u_-) + \ab(u_-,u') \, (u_- - u') - \ab(u_+,u') \, (u_+ - u')
\\
& = f(u_+) - f(u_-) + (f(u_-) - f(u')) - (f(u_+) - f(u')) = 0, 
\endaligned 
$$
we obtain the identity 
\be
\label{scalar10}
\bigl(   \ab(u_-, u_+) - \ab(u_-, u') \bigr) \,  (u_- - u')
= - 
\bigl( \ab(u_-, u_+)  - \ab(u_+, u') \bigr) \,  (u_+ - u'). 
\ee

Let us recall Definition~\ref{def21} which provides the classification of shocks. 
The identity \eqref{scalar10}
shows that the shock in the function $\abf$ is a compressive or rarefaction shock 
if $u_\pm - u'$ have opposite signs; it is undercompressive if $u_\pm - u'$ have the same sign.
(In fact, \eqref{scalar10} is nothing but the identity \eqref{2.13} already  derived for 
general systems, but for scalar equations \eqref{2.13} reduces to saying that $\gamma_1^+ = \gamma_1^-$.) 

Next, it is interesting to observe that the identity \eqref{scalar10} can be cast in the new form
\be
\label{scalar11}
\bigl( \ab(u_+, u') - \ab(u_-, u') \bigr) \,  (u_- - u')
= 
\bigl( \ab(u_+, u') - \ab(u_-, u_+) \bigr) \,  (u_+ - u_-). 
\ee
This second identity allows us now to identify the nature of undercompressive shocks
(either slow or fast) when $u' - u_\pm$ have the same sign, and the nature of 
compressive and rarefaction shocks when $u' - u_\pm$ have opposite signs. 
\end{proof}

A result similar to Lemma~\ref{Lemma3.2}, including the key identities \eqref{scalar10} and \eqref{scalar11}, 
will be established for fluid dynamics equations 
in the course of the proof of the following theorem
and, later, a version valid for general systems will be established as well. 

\begin{theorem}[Averaged matrix of scalar and fluid dynamics equations]
\label{Prop3.1} 
Let~$u=u(t,x)$ and $u'=u'(t,x)$ be solutions with bounded variation and {\rm arbitrary large amplitude} to 
\begin{enumerate}
\item a scalar conservation law with general flux-function,  or  
\item the Lagrangian $p$-system of isentropic fluid dynamics with general volume-pressure law 
$p=p(v)$ satisfying solely the hyperbolicity condition $\del_v p(v)<0$.
\end{enumerate}
Then an entropy admissible shock of $u$ or $u'$ can not be a rarefaction-shock of the averaged matrix $\Abf(t,x) := \Ab(u,u')(t,x)$. 
\end{theorem} 
 
\begin{remark}
1. As observed earlier, the averaged matrix is not unique in general (except when $N=1$). 
For the $p$-system we use the averaged matrix defined by the general formula \eqref{matrix1}. 
Later we will see (e.g.~in Section~\ref{Section6}) that for the Euler equations
the averaged matrix must be based on the special structure of the equations. 
 
2. Theorem~\ref{Prop3.1} provides us with a key property required for $L^1$ stability. 
Of course, approximate solutions $u,u'$ (for instance constructed by front tracking) generally admit 
non-admissible shocks; however, the strength of such shocks vanish 
in the limit (as the approximation parameters tend to zero) so that they 
do not prevent the $L^1$ stability estimate to hold in the limit.   
\end{remark}

\begin{proof} The statement for scalar conservation laws is consequence of Lemma~\ref{Lemma3.2}
and the entropy condition. 
Namely, consider a point $(t,x)$ where the solution admits an entropy admissible  
shock $(u_-,u_+)$ while $u'$ is locally constant. 
Recall Oleinik's entropy inequalities \cite{Oleinik} 
\be
{f(v) - f(u_-) \over v - u_-} \geq {f(u_+) - f(u_-) \over u_+ - u_-} 
\quad 
\text{ for all } v \text{ between $u_-$ and $u_+$.} 
\label{3.8}
\ee 
According to Lemma~\ref{Lemma3.2}, to exclude rarefaction-shocks we must exclude the first case in \eqref{3.4}
 where $u' - u_\pm$ have opposite signs and, simultaneously, 
\be
(u' - u_-) \, (\abf_+ - \abf_-) \, (u_+ - u_-) < 0.  
\label{3.9}
\ee
Suppose \eqref{3.9} hold with, for definiteness, $u_+ < u_-$. Hence $u_+ < u' < u_-$, so 
$u_+ - u_- < 0$ and $u' - u_- < 0$. It follows from \eqref{3.9} that $\abf_+ - \abf_- > 0$.
On the other hand, according to Oleinik's condition \eqref{3.8}
$$
\abf_- \geq \abf_+, 
$$ 
which leads to a contradiction. Consequently, the first case 
in \eqref{3.4} is impossible, and the shock is compressible or undercompressive, but can not be 
a rarefaction-shock. This completes the proof of the theorem for a scalar conservation law. 

Consider next the $p$-system
\be
\label{psystem}
\del_t u + \del_x p(v) = 0, \qquad \del_t v - \del_x u = 0, 
\ee
where $u$ and $v$ represent the fluid velocity and the (positive) specific volume, respectively. 
The pressure function $p=p(v)$ is assumed to be decreasing, which ensures that \eqref{psystem}
is strictly hyperbolic. 
Given two solutions $(u,v)$, $(u',v')$ of the $p$-system and setting $\varphi := u - u'$, $\psi := v - v'$, we obtain 
$$
\del_t \varphi - \del_x \big( \cb(v,v')^2 \, \psi) = 0, 
\qquad 
\del_t \varphi - \del_x \psi = 0, 
$$ 
where 
$$
\cb(v,v') = \sqrt{- {p(v) - p(v') \over v - v'}}.  
$$
Therefore, the averaged matrix \eqref{matrix1} associated with the $p$-system is   
\be
\label{matrixp}
\Ab(u,v; u',v') = \Ab(v, v') := \begin{pmatrix} 0 & - \cb(v,v')^2
\\
- 1   &  0
\end{pmatrix}, 
\ee
where $\lamb_j(u,v; u,v'):= \pm\cb(v,v')$ are the averaged wave speeds and correspond to 
the physical wave speeds $\lambda_j(u,v) = \pm c(v):=\pm\sqrt{-\del_v p(v)}$. 

Observe that the averaged speed of the $p$-system only depends upon a single variable, i.e.~the 
specific volume $v,v'$. Considering a $2$-shock wave connecting $(u_-,v_-)$ to $(u_+,v_+)$
at the (positive) speed $\cb(v_-,v_+)$ and denoting by $(u',v')$ an arbitrary constant, we find the identity 
\be
\label{p43}
\aligned 
& \bigl( \cb(v_-, v_+) - \cb(v_-, v') \bigr) \,  (v' - v_- )
\\
& = \kappa(v_-,v_+,u') \bigl( \cb(v_-, v_+) - \cb(v_+, v') \bigr) \,  (v' - v_+ ), 
\endaligned 
\ee
where 
$$
\kappa(v_-,v_+,v') := \bigl( \cb(v_-, v_+) + \cb(v_+, v') \bigr) \,  \bigl( \cb(v_-, v_+) + \cb(v_-, v') \bigr)^{-1}. 
$$
This shows that the $2$-shock of the averaged matrix 
\be
\Abf(t,x) = \begin{cases} 
\Ab(v_-,v'),  &  x <  \cb(v_-, v_+) \, t, 
\\
\Ab(v_+,v'),  &  x >  \cb(v_-, v_+) \, t, 
\end{cases}
\label{3.299}
\ee 
is a compressive or rarefaction shock if $v_\pm- v'$ have opposite signs, 
while it is undercompressive if $v_\pm - v'$ have the same sign.

We can also derive the second identity 
\be
\label{p44}
\aligned 
& \bigl( \cb(v_+, v') - \cb(v_-, v') \bigr) \,  (v_- - v')
\\
& = \kappa(v',v_+,v_-) \, \bigl( \cb(v_+, v') - \cb(v_-, v_+) \bigr) \,  (v_+ - v_-),
\endaligned 
\ee
which allows us to distinguish between compressive and rarefaction $2$-shocks, on one hand, and, on the other hand, 
between 
slow or fast undercompressive $2$-shocks. This leads us to an analogue of Lemma~\ref{Lemma3.2} for the $p$-system: 
setting 
$$
\rho(v_-,v_+,v') :=  (v_- - v') \, (v_+ - v_-), 
$$
we have : 
\begin{enumerate} 
\item If $v_\pm - v'$ have the {\rm same sign}, 
\be 
\rho(v_-,v_+,v') \, (\cb(v_+,v') - \cb(v_-,v') 
\begin{cases} 
\leq 0, & \Scal_2^\Abf,     
\\ 
\geq 0, & \Fcal_2^\Abf. 
\end{cases}  
\label{3.399}  
\ee
\item If $v_\pm - v'$ have {\rm opposite signs}, 
\be
\rho(v_-,v_+,v') \, (\cb(v_+,v') - \cb(v_-,v')) 
\begin{cases} 
\leq 0, & \Rcal_2^\Abf,     
\\ 
\geq 0, & \Lcal_2^\Abf. 
\end{cases}  
\label{3.499}  
\ee 
\end{enumerate} 
In turn, recalling that Wendroff's entropy criterion for $2$-shocks of the $p$-system \cite{Wendroff} 
\be
\cb(v,v_-) \geq \cb(v_-,v_+) 
\quad 
\text{ for all } v \text{ between $v_-$ and $v_+$,} 
\label{Wend}
\ee 
we can check that that the first case in \eqref{3.499} can not arise if the shock is entropy admissible. 
We conclude that, for the $p$-system,
admissible $2$-shocks can not be rarefaction-shocks of the averaged matrix \eqref{matrixp}. 
A similar conclusion holds for shock waves in the first characteristic family, and this completes the proof of the theorem. 
\end{proof}

\subsection{Genuinely nonlinear systems}
 
In the case of genuinely nonlinear characteristic fields we can derive a monotonicity property 
of the eigenvalues of the averaged matrix. 

\begin{lemma}[Shock speed monotonicity property] 
\label{Lemma3.3} 
Let \eqref{1.1} be a strictly hyperbolic system of conservation laws with flux $f:\Bone \to \RN$ 
admitting a genuinely nonlinear eigenvalue $\lambda_i$.  Then there exists $\delta_1 < \delta_0$
such that the following property holds for every averaged matrix satisfying \eqref{matrix0}. 
For arbitrary $u_-, u' \in \Bone$ the averaged wave speed $\lamb_i(\cdot,u')$ is a strictly 
monotone function along the $i$-shock curve issuing from $u_-$.  
More precisely, if an $i$-shock $(u_-,u_+)$ satisfies Lax shock inequalities 
\be 
\lam_i(u_-) >  \lamb_i(u_-, u_+) > \lam_i(u_+), 
\label{3.5}
\ee
then the averaged speed satisfies the inequality 
\be
\lamb_i(u_-,u') > \lamb_i(u_+,u'). 
\label{3.6}
\ee
Furthermore, this holds for the $p$-system of Lagrangian fluid dynamics 
even for shocks with {\rm arbitrary large amplitude,}
provided the equation of state is genuinely nonlinear (that is, $\del_{vv} p(v) > 0$). 
\end{lemma}

\begin{proof} The right-hand state $u_+=u_+(\eps)$ can be viewed as a function of $u_-$ and a parameter $\eps$
varying in the neighborhood of $0$, with 
$$ 
u_+(\eps) = u_- + \eps \, r_i(u_-) + O(\eps^2). 
$$ 
Then, we compute 
$$
\aligned 
\lamb_i(u_+,u') - \lamb_i(u_-,u')
& = \eps \, \nabla_1 \lamb_i(u_-, u') \cdot r_i(u_-) + O(\eps^2)
\\
& = \eps \, \nabla_1\lamb_i(u_-,u_-) \cdot r_i(u_-) + O(\eps^2) + O(\eps \, |u' - u_-|) > 0. 
\endaligned  
$$
Using now that $\lamb_i(u,u) = \lam_i(u)$ and the symmetric property $\Ab(u,u') = \Ab(u',u)$, 
we obtain $2 \, \nabla_1 \lamb_i \cdot r_i = \nabla \lam_i \cdot r_i >0$.  
Provided $|\eps| + |u'- u_-| \lesssim \delta_1$ is sufficiently small, 
we conclude that the function $\lamb_i(\cdot,u')$ is strictly monotone along the shock curve. 

If the shock $(u_-,u_+)$ satisfies the entropy inequalities \eqref{3.5} 
and the normalization $\nabla\lam_i \cdot r_i >0$ is chosen, then $\eps < 0$ 
and we conclude from the above calculation that the averaged speed decreases from $u_-$ to $u_+$. 

For the $p$-system, using the notation introduced earlier we can write
$$
{\del \cb \over dv} (v,v') = - {1 \over \cb(v,v')} \int_0^1 \del_{vv}p(v' + \theta(v-v')) \theta \, d\theta, 
$$
which keeps a constant sign if $\del_{vv} p$ has a constant sign. Therefore, under the genuine nonlinearity assumption, 
the averaged speed of the $p$-system is strictly monotone along each Hugoniot curve. 
\end{proof} 

We deduce that: 

\begin{theorem}[Averaged matrix of genuinely nonlinear systems]
\label{GNLmonotone}
Let $u=u(t,x)$ and $u'=u'(t,x)$ be (small amplitude) solutions with bounded variation to 
a system of conservation laws \eqref{1.1} with genuinely nonlinear flux $f$, 
and let $\Ab$ be an averaged matrix satisfying \eqref{matrix0}. 
Then, an entropy admissible shock of $u$ or $u'$ can not be a rarefaction-shock 
of the averaged matrix $A(t,x) := \Ab(u, u')(t,x)$. 
\end{theorem}

For scalar conservation laws and fluid dynamics equations, 
the monotonicity property exhibited in Lemma~\ref{3.3} provides us with a 
second proof  of Theorem~\ref{Prop3.1} which, however, requires genuine nonlinearity. 

\begin{proof} We rely on Lemma~\ref{3.3} and consider a point $(t,x)$ of jump of the BV function 
$u$: the limit $u_-$ is connected to $u_+$ by an $i$-shock wave with speed 
$\lamb_i(u_-, u_+)$ satisfying Lax shock inequalities \eqref{3.5}. If the $BV$ solution $u'$ is continuous at  
$(t,x)$ then the inequalities characterizing a rarefaction-shock, 
$$
\lamb_i(u_-, u') < \lamb_i(u_-, u_+) < \lamb_i(u_+, u'), 
$$ 
clearly contradict the monotonicity property \eqref{3.6}. 
\end{proof}


\subsection{Systems with general flux} 
\label{subsection3.3}

For characteristic fields that are not genuinely nonlinear, 
a suitable version of Lemma~\ref{Lemma3.3} and Theorem~\ref{GNLmonotone} will now be established, 
by restricting attention to ``robust'' patterns $(u_-,u_+,u')$,
corresponding to ``strongly dominant'' characteristic components, in a sense 
introduced below. 

We begin with a preliminary result. 
 
\begin{lemma}
\label{Lemma3.4} 
Given  an $i$-shock wave connecting the states $u_-, u_+ \in \Bzero$ and given $u' \in \Bzero$, 
define 
$$
\aligned 
& \lamb_j^\pm = \lamb_j(u_\pm,u'), \quad \lamb = \lamb_i(u_-,u_+), \quad 
\rb_j^\pm: = \rb_j(u_\pm,u'),
\\
& \lb_j^\pm: = \lb_j(u_\pm,u'), \quad \eps^\pm =  l_i^\pm \cdot (u_+ - u_-), 
\\
& u_\pm - u' =: \sum_j \alpha_j^\pm \, \rb_j^\pm. 
\endaligned 
$$ 
Then, the folllowing identities hold
\be
\aligned 
& \bigl( \lamb_i^+ - \lamb \bigr) \, \eps^- = - \bigl( \lamb_i^+ - \lamb_i^- \bigr) \alpha_i^- + O(1) \, \Omega,  
\\
& \bigl( \lamb - \lamb_i^- \bigr) \, \eps^+ = \bigl( \lamb_i^+ - \lamb_i^- \bigr) \alpha_i^+ + O(1) \, \Omega, 
\endaligned 
\label{3.10}
\ee
where the remainder $\Omega$ satisfies 
$$
\aligned 
\Omega 
& = O(1) \, | \rb_i^+ - \rb_i^-| \, |\lamb_i^- - \lamb | \, |\alpha_i^-| + O(\eps^-) \, \sum_{k\neq i} |\alpha_k^-|. 
 \endaligned  
$$
\end{lemma}

In the course of the proof given below, we will derive the more general identities 
\be
\label{relat33}
\aligned 
& (\lamb_i^+  - \lamb_i^- ) \, (u_+ - u') - (\lamb  - \lamb_i^-) \, (u_+ - u_-) 
\\
& = ( \lamb_i^+ - \lamb_i^- ) \, (u_- - u') + ( \lamb_i^+ - \lamb ) \, (u_+ - u_-) 
= O(1) \,  \Omega, 
\endaligned 
\ee 
together with the following expression for the remainder
\be
\label{relat34}
\aligned 
\Omega: 
& =  \sum_{j \neq i}|\bigl( \lamb - \lamb_j^+ \bigr)  \, \alpha_j^+ - \bigl( \lamb - \lamb_j^- \bigr)  \, \alpha_j^-| 
      + O(\eps^-) \sum_{j \neq i} |\alpha_j^-|.  
 \endaligned  
\ee

\begin{proof}  We can write 
\be
\label{identity5}
\aligned 
& \bigl(  \lamb_i(u_+,u')  - \lamb_i(u_-,u') \bigr) \, (u_- - u') +  \bigl(  \lamb_i(u_+,u')  - \lamb_i(u_-,u_+) \bigr) \, (u_+ - u_-)
\\
& = \bigl(  \lamb_i(u_+,u')  - \lamb_i(u_-,u') \bigr) \, (u_+ - u')  - \bigl(  \lamb_i(u_-,u_+)  - \lamb_i(u_-,u') \bigr) \, (u_+ - u_-) 
\\
& = \Omegat,
\endaligned 
\ee
where, using that $u_+ $ lies on the Hugoniot curve from $u_-$,  
$$
\aligned 
\Omegat := & - \big( f(u_+) - f(u') - \lamb_i(u_+,u') \, (u_+ - u') \big) 
                      +  \big( f(u_-) - f(u') - \lamb_i(u_-,u') \, (u_- - u') \big)
\\
= &  - \big( \Ab(u_+,u') - \lamb_i(u_+,u') \big) \, (u_+ - u')  
       + \big( \Ab(u_-,u') - \lamb_i(u_-,u') \big) \, (u_- - u')
\\
= & - \sum_{j \neq i} \big( \lamb_i(u_+,u') - \lamb_j(u_+,u') \big) \, \alpha_j(u_+, u') \, \rb_j(u_+,u')
\\ 
   & + \sum_{j \neq i} \big( \lamb_j(u_-,u')  - \lamb_j(u_-,u') \big) \, \alpha_j(u_-, u') \, \rb_j(u_-,u'). 
\endaligned 
$$ 
Therefore, introducing the notation 
$$
\gamma_j^\pm = \gamma_j^\pm(u_-,u_+,u') := \bigl( \lamb - \lamb_j^\pm \bigr)  \, \alpha_j^\pm,   
$$
the remainder takes the form 
$$
\aligned
\Omegat = & - \sum_{j \neq i} {\lamb_j^+ - \lamb_i^+ \over \lamb - \lamb_j^+} \, \gamma_j^+ \, \rb_j^+
                      +  \sum_{j \neq i} {\lamb_j^- - \lamb_i^- \over \lamb - \lamb_j^-} \, \gamma_j^- \, \rb_j^-.   
\endaligned 
$$
Observing the decomposition 
$$
\aligned
\Omegat & = - \sum_{j \neq i} {\lamb_j^+ - \lamb_i^+ \over \lamb - \lamb_j^+} \, (\gamma_j^+ - \gamma_j^-) \, \rb_j^+
                     -  \sum_{j \neq i} \Big(  {\lamb_j^+ - \lamb_i^+ \over \lamb - \lamb_j^+} \, r_j^+  - 
         {\lam_j^- - \lamb_i^- \over \lamb - \lamb_j^-} \, r_j^- \Big) \,  \gamma_j^-
         \\
         & = O(1) \sum_{j \neq i}|\gamma_j^+ - \gamma_j^-| + O(\epsb^\Abf) \sum_{j \neq i} |\gamma_j^-|, 
\endaligned 
$$
we arrive at the general statement \eqref{relat33}-\eqref{relat34}. 

Multiplying now \eqref{identity5} by $l_i^-$ and $l_i^+$ respectively, we obtain 
$$
\aligned 
& (\lamb_i^-  - \lamb ) \, l_i ^+ \cdot (u_+ - u_-) - (\lamb_i^+  - \lamb_i^-) \, \alpha_i^+ = - l_i ^+ \cdot \Omegat, 
\\
& (\lamb_i^+  - \lamb) \, l_i ^- \cdot (u_+ - u_-) - ( \lamb_i^+ - \lamb_i^- ) \, \alpha_i^- = - l_i ^- \cdot \Omegat,
\endaligned 
$$ 
which gives \eqref{3.10}. Next, using $u_\pm - u' =: \sum_j \alpha_j^\pm \, \rb_j^\pm$ we find
$$
- ( \lamb_i^+  - \lamb) \, \sum_j \alpha_j^+ \, \rb_j^+ 
+ ( \lamb - \lamb_i^- ) \, \sum_j \alpha_j^+ \, \rb_j^+ = \Omegat,
$$
thus 
$$
\gamma_i^+ \, \rb_i^+ - \gamma_i^- \, \rb_i^- 
= \Omegat - (\lamb_i^+  - \lamb) \, \sum_{j \neq i} \alpha_j^+ \, \rb_j^+ 
   + (  \lamb - \lamb_i^- ) \, \sum_{j \neq i} \alpha_j^- \, \rb_j^-,  
$$  
which leads to the desired estimate \eqref{3.10} along the same lines as in the proof of Lemma~\ref{jump}. 
\end{proof}

The relations \eqref{3.10} are the analogue, to the case of systems, of the identity \eqref{scalar11} found earlier for scalar equations. 
They allow us to relate the signs of $\lamb_i^\pm - \lamb$ with the sign of $\lamb_i^+ - \lamb_i^-$. 
This is so only when the reminder $\Omega$ can be neglected, and this motivates us to propose 
the following definition.

\begin{definition}  
\label{strong}
With the notation in Lemma~\ref{Lemma3.4}, 
the $i$-characteristic component  across an $i$-shock is said to be {\rm strongly dominant} if 
its left- and right-hand limits $\alpha_i^\pm$ satisfy 
$$ 
\kappa_2 \, | \lamb_i^+ - \lamb_i^- | \, | \alpha_i^\pm| \geq
 | \rb_i^+ - \rb_i^-| \, |\lamb_i^- - \lamb | \, |\alpha_i^-| + |u_+ - u_-| \, \sum_{k\neq i} |\alpha_k^-|.  
$$
In this case, the shock $(u_-,u_+)$ is said to be {\rm robust relatively to the state $u'$.}  
\end{definition}

Following \cite{IguchiLeFloch} we parametrize the $i$-wave curves with a global parameter $\mu_i(u)$ 
such that 
$$
\nabla \mu_i \cdot r_i(u) > 0. 
$$ 
The hypersurfaces along which $\mu_i$ is constant determine a smooth foliation of the phase space
and all $i$-curves are transverse to them. Observe that 
$$
\aligned 
\sgn \big( \lb_i(u_-,u') \cdot (u_+ - u_-) \big) 
& = \sgn \big( \lb_i(u_+,u') \cdot (u_+ - u_-) \big) 
\\
& = \sgn \big( \mu_i(u_+) - \mu_i(u_-) \big),   
\endaligned 
$$
when $u_+$ lies on the Hugoniot curve from $u_-$ since this curve is tangent to the $i$-eigenvector at $u_-$. 

From the identities in Lemma~\ref{Lemma3.4} we deduce immediately: 

\begin{lemma}[Characterization of the nature of robust shocks]  
\label{Prop3.5} 
With the notation in Lemma~\ref{Lemma3.4}, the nature of a robust $i$-shock of the averaged matrix 
$$
\Abf(t,x) = \begin{cases} 
\Ab(u_-,u'), & x < \lamb_i(u_-, u_+) \, t, 
\\
\Ab(u_+,u'), & x > \lamb_i(u_-, u_+) \, t,
\end{cases} 
$$ 
is determined by the signs of 
$\rho_i(u_-,u_+,u') := \alpha_i(u_-,u')  \, (\mu_i(u_+) - \mu_i(u_-))$
and the jump $(\lamb_i(u_+,u') - \lamb_i(u_-,u') )$ of the $i$-eigenvalue, as follows: 
\begin{enumerate}
\item If $\alpha_i(u_\pm,u')$ have the same sign, 
\be 
\rho_i(u_-,u_+,u') \, (\lamb_i(u_+,u') - \lamb_i(u_-,u') )
\begin{cases} 
< 0,     &  
 \Scal_i^\Abf \text{ and robust,} 
\\ 
\\ 
> 0,     & \Fcal_i^\Abf \text{ and robust.}  
\end{cases} 
\label{3.13}
\ee
\item If $\alpha_i(u_\pm,u')$ have opposite signs, 
\be
\rho_i(u_-,u_+,u') \, (\lamb_i(u_+,u') - \lamb_i(u_-,u') )
\begin{cases} 
> 0,     & \Lcal_i^\Abf \text{ and robust,}
\\ 
\\ 
< 0,     & \Rcal_i^\Abf \text{ and robust.} 
\end{cases} 
\label{3.14}
\ee
\end{enumerate} 
\end{lemma}

\begin{remark} For scalar conservation laws every shock is robust (since $\Omega \equiv 0$) 
and we can recover Lemma~\ref{Lemma3.2} from Lemma~\ref{Prop3.5}. 
\end{remark}

Using the above result together with the entropy condition which for general system
imposes that $(u_-,u_+)$ is entropy admissible \cite{Liu1,Liu2} if and only if the shock speed 
$\lamb_i(u_-,\cdot)$ achieves at the point $u_+$ its minimal value along the Hugoniot curve, 
that is 
\be
\lamb_i(u_-, u_+) \leq \lamb_i(u_-,v)
\label{3.111} 
\ee
for every state $v$ along the Hugoniot curve between $u_-$ to $u_+$. Combining now the conclusions of Lemma~\ref{Prop3.5}
with the entropy inequality \eqref{3.111} we arrive at:

\begin{theorem}[Averaged matrix for systems with general flux]
\label{general}
Let $u=u(t,x)$ and $u'=u'(t,x)$ be (small amplitude) solutions with bounded variation to the strictly hyperbolic 
system of conservation laws \eqref{1.1} with general flux $f$. 
Then, a robust and entropy admissible shock of $u$ or $u'$ can not be a rarefaction-shock of the averaged matrix 
$A(t,x) := \Ab(u, u')(t,x)$. 
\end{theorem}


\section{Stability for a class of linear hyperbolic systems}
\label{section4} 

We are now in a position to derive an $L^1$ stability property for a large class of systems, 
by assuming the existence of a weight $\wbf = \wbf(t,x)$ satisfying certain constraints 
that we specify. Consider a solution $\psi=\psi(t,x)$ of a uniformly hyperbolic system with bounded variation \eqref{1.5}, 
together with its characteristic components $\alpha_j, \beta_j$ defined earlier in Section~\ref{LI-0}.  
Relying on Lemma~\ref{weight} on the evolution of the weighted norm, 
our goal is to establish that the right-hand side of \eqref{2.8} is non-positive. 
At an $i$-shock, Lemma~\ref{jump} and, more precisely, \eqref{beta} 
shows that the components $|\beta_j^-|$ and $|\beta_j^+|$ coincide 
``up to first-order'', with error terms proportional to $| \Abf^+ - \Abf^- |$ or $|r^\Abf_{i+} - r^\Abf_{i-}|$. 
The weight $\wbf$ should have jumps that precisely compensate the effect of these error terms. 
In the present section,  we assume that the matrix-field $\Abf$ satisfies the condition 
\be
\label{eigenv}
| r_{i+}^\Abf -  r_{i-}^\Abf | \lesssim | \lam_{i+}^\Abf - \lam_{i-}^\Abf | 
\ee
at every jump point, and we formulate the jump conditions in term of the eigenvalues $(\lam^\Abf_{i+} - \lam^\Abf_{i-})$
rather than the eigenvectors. 

\begin{definition}  
\label{strong2}
With the above notation, the $i$-characteristic component  across an $i$-shock is said to be {\rm strongly dominant} if 
its left- and right-hand limits $\alpha_i^\pm$ satisfy 
$$ 
\kappa_2 \, | \lamb_i^+ - \lamb_i^- | \, | \alpha_i^\pm| \geq | \Abf^+ - \Abf^- | \, \sum_{k\neq i} |\alpha_k^-|.  
$$ 
\end{definition}

We seek for conditions on $\wbf_j^+ - \wbf_j^-$ ensuring that the right-hand side of \eqref{2.8} is non-positive.
In view of \eqref{beta}, error terms associated with 
non-dominant components can be estimated by the dominant ones 
\be
\label{2.rough}
\sum_{j} |\beta_j^-| \leq C \, \sum_{j \text{ dominant}} |\beta_j^-|, 
\ee 
and no condition will be imposed on the jump $\wbf_j^+ - \wbf_j^-$ of non-dominant components. 
On the other hand, in view of the sign properties \eqref{sign-j}-\eqref{sign-i} of the characteristic flux $\beta_j^\pm$, 
it is natural to require that the weight $\wbf$ satisfies \eqref{updown}. 

More precisely, we impose that at each $(t,x) \in \Jcal_i^\Abf$ and for $j \neq i$,  
\be
\label{contraintes-j}
\wbf_j^+ - \wbf_j^- 
\begin{cases}  
\geq \hskip.3cm K \, | \Abf^+ - \Abf^- |,  & j<i \text{ and $j$ dominant,} 
\\ 
\leq - K \, | \Abf^+ - \Abf^- |,  & j>i \text{ and $j$ dominant}  
\end{cases}
\ee  
for a sufficiently large constant $K>0$. 

Specifying the jump $\wbf_i^+ - \wbf_i^-$ of the $i$-component across an $i$-shock is more delicate. 
Observe that no condition is necessary on $\Lcal_i^\Abf$ (compressive shocks) since both flux 
$\beta_i^\pm$ have a favorable
sign. On the other hand, it is hopeless to try to impose a condition on $\Rcal_i^\Abf$ (rarefaction shocks) since both $\beta_i^\pm$
have an unfavorable sign. Recall here Definition~\ref{strong} and Lemma~\ref{Prop3.5} 
which show that only strongly dominant rarefaction-shocks have been characterized, however. This motivates us to 
allow only such waves in the right-hand side of the key estimate \eqref{4.66} below. 
 
Next, we focus on undercompressive shocks associated with dominant components. 
At this juncture, we recall Definition~\ref{strong} and Lemma~\ref{Prop3.5} 
which provide us with the nature (slow or fast) of an undercompressive shock, but
under the assumption that the characteristic components are strongly dominant. 
This motivates the following constrains on the weight:   
\be
\label{contraintes-i1}
\wbf_i^+ - \wbf_i^- 
\begin{cases}  
\leq - K \,  | \lam_{i+}^\Abf -  \lam_{i-}^\Abf |,  & \Scal_i^\Abf \text{ and strongly dominant,} 
\\ 
\geq  \hskip.3cm  K \, | \lam_{i+}^\Abf -  \lam_{i-}^\Abf |, & \Fcal_i^\Abf \text{ and strongly dominant,}
\end{cases}
\ee 
while for dominant but not strongly dominant component we require the weaker condition 
\be
\label{contraintes-i2}
\wbf_i^+ - \wbf_i^- 
\begin{cases}  
\leq  \hskip.3cm  K \,  | \lam_{i+}^\Abf - \lam_{i-}^\Abf |, 
& \Scal_i^\Abf \text{ and $i$ dominant,} 
\\ 
\geq  - K \, | \lam_{i+}^\Abf - \lam_{i-}^\Abf |, & \Fcal_i^\Abf \text{ and $i$ dominant.} 
\end{cases}
\ee 
No condition is required on $\wbf_i$ if the $i$-characteristic component is not dominant. 

\begin{proposition} 
\label{LI-main} 
Let $\Abf = \Abf(t,x)$ be a uniformly hyperbolic matrix with bounded variation 
satisfying \eqref{eigenv}, 
Then, given any piecewise constant solution  $\psi=\psi(t,x)$ of \eqref{1.5}
and a piecewise constant weight $\wbf = \wbf(t,x)$ satisfying the general conditions stated in Section~\ref{LI-0}, together with
 \eqref{contraintes-j}--\eqref{contraintes-i2}, the associated weighted norm satisfies 
$$
{d \over dt}  \| \psi(t)\|_{\wbf(t)} \lesssim \sup_{x \in \Kcal_i^\Abf(t)} | \psi^+(t,x) - \psi^-(t,x)|, 
$$ 
where $\Kcal_i^\Abf \subset \Rcal_i^\Abf$ denote all $i$-rarefaction-shocks associated with strongly dominant 
characteristic components.  Consequently, the solution $\psi$ satisfies the $L^1$ stability estimate
\be
\| \psi(t) \|_{L^1(\R)}  
\lesssim  \|\psi(0)\|_{L^1(\R)} + \int_0^t \sup_{x \in \Kcal_i^\Abf(s)} | \psi^+(s,x) - \psi^-(s,x)| \, ds, \qquad t \in \R_+. 
\label{4.66}
\ee 
\end{proposition}

We conclude that the $L^1$ stability property \eqref{1.6} holds when $\Abf$ does not admit 
strongly dominant rarefaction shocks. 
This result will be used in Section~\ref{Section6} in combination with our earlier conclusion (Theorem~\ref{general}) 
that the averaged matrix of two entropy solutions of a system of conservation laws does not admit 
strongly dominant rarefaction-shocks.   

It is worth pointing out also that Proposition~\ref{LI-main} 
remains valid for {\sl approximate} solutions that satisfy the equations \eqref{1.5} up to a measure source-term, 
provided the total mass of the source is added to the right-hand side of \eqref{4.66}.  
See \cite{LeFloch-book} for details.

\begin{proof} We use Lemma~\ref{weight} and estimate each term in the right-hand side of \eqref{2.8}: 
$$ 
B:= \sum_j \beta_j^- \, \wbf_j^- + \beta_j^+ \, \wbf_j^+ 
$$ 
associated with a given $i$-shock. In the following, we will often make use of \eqref{2.rough}. 
Using the expansion \eqref{beta} we can write 
$$
\aligned 
B = 
& \sum_j \bigl( \beta_j^- \, \wbf_j^- + \sgn(\beta_j^+) \, |\beta_j^-|  \, \wbf_j^+ \bigr)  
   \\
   &  + O(| r_{i+}^\Abf - r_{i-}^\Abf |) |\beta_i^-|   + O(| \Abf^+ - \Abf^- |) \sum_{k\neq i} |\beta_k^-|,  
\endaligned 
$$
so that by \eqref{sign-j} and \eqref{eigenv} 
$$ 
\aligned 
B = 
& \bigl( \wbf_i^- \, \sgn(\beta_i^-) + 
          \wbf_i^+ \, \sgn(\beta_i^+) \bigr) \,  |\beta_i^-| - \sum_{j < i}  \bigl( \wbf_j^+ - \wbf_j^- \bigr) \, |\beta_j^-|   
\\
& 
  - \sum_{j > i}   \bigl( \wbf_j^- - \wbf_j^+ \bigr) \, |\beta_j^-|      + O(| \lam_{i+}^\Abf - \lam_{i-}^\Abf |) |\beta_i^-| 
   + O(| \Abf^+ - \Abf^- |) \sum_{k\neq i} |\beta_k^-|. 
\endaligned 
$$

Thanks to \eqref{contraintes-j}, for $j \neq i$ the dominant $j$-components lead 
to a favorable sign of the jump ${\wbf_j^+ - \wbf_j^-}$, while non-dominant components for $j \neq i$ 
can be collected in the remainder using \eqref{2.rough}. We obtain 
$$ 
\aligned 
B \leq & \bigl( \wbf_i^- \, \sgn(\beta_i^-) 
           + \wbf_i^+ \, \sgn(\beta_i^+ ) \bigr) \, |\beta_i^-| 
      \\
      &   - K \, | \Abf^+ - \Abf^- | \hskip-.2cm \sum_{k \neq i \atop \text{dominant}}  \hskip-.2cm |\beta_k^-|    
        + O(| \lam_{i+}^\Abf - \lam_{i-}^\Abf |) |\beta_i^-|   + O(| \Abf^+ - \Abf^- |) \sum_{k\neq i} |\beta_k^-|. 
\endaligned 
$$
The sum over dominant components allows us
to suppress the remainder, provided that $K$ is sufficiently large (so that $K \, | \Abf^+ - \Abf^- |$ dominates $O(| \Abf^+ - \Abf^- |)$). 
We arrive at the inequality 
\be
\label{2.31}
\aligned 
B \leq  
& \bigl( \wbf_i^- \, \sgn(\beta_i^-)  + \wbf_i^+ \, \sgn(\beta_i^+) \bigr) \, |\beta_i^-|
\\
& + O(| \lam_{i+}^\Abf - \lam_{i-}^\Abf |) \, |\beta_i^-| - {K \over 2} \, | \Abf^+ - \Abf^- | \, \sum_{k \neq i} |\beta_k^-|. 
\endaligned 
\ee 
It remains to deal with the term $|\beta_i^-|$, which can be assumed to be dominant, since 
non-dominant components $|\beta_i^-|$ can be handled with \eqref{2.rough}. 
We will now distinguish between the cases of compressive, undercompressive, and rarefaction shocks. 

\vskip.3cm 

\noindent{\bf Case 1.} If the $i$-shock is compressive, then we have $\sgn(\beta_j^-) 
=  \sgn(\beta_j^+) = -1$ and therefore
$$ 
B \leq - 2 \, w^{\min} \, |\beta_i^-| + O(| \lam_{i+}^\Abf - \lam_{i-}^\Abf |) \,|\beta_i^-| - { K \over 2} \, | \Abf^+ - \Abf^- | \, \sum_{j \neq i} |\beta_j^-|.
$$
So we obtain 
\be
\label{2.compress}
B \leq -  w^{\min} \, |\beta_i^-| - {K \over 2} \, | \Abf^+ - \Abf^- | \, \sum_{j \neq i} |\beta_j^-| 
\qquad \text{ in }  \Lcal_i^\Abf.  
\ee 

\vskip.3cm 
 
\noindent{\bf Case 2. } If the $i$-shock is undercompressive, we will prove that  
\be
B \leq - {K \over 3} \, | \lam_{i+}^\Abf - \lam_{i-}^\Abf | \, |\beta_i^-|  - {K \over 3} \, | \Abf^+ - \Abf^- | \, \sum_j |\beta_j^-| 
\qquad \text{ in }  \Scal_i^\Abf \cup \Fcal_i^\Abf.  
\label{2.33}
\ee
Here, we have $\sgn(\beta_i^- ) = - \sgn(\beta_i^+)$.

First of all, if the $i$-component is strongly dominant, then by \eqref{contraintes-i1} 
we have determined the jump $\wbf_i^+ - \wbf_i^-$ in such a way that
it can compensate the error term $| \lam_{i+}^\Abf - \lam_{i-}^\Abf | \, |\beta_i^-| $, provided 
the constant $K$ is chosen to be sufficiently large 
so that $K \, | \Abf^+ - \Abf^- |$ is larger than $O(| \Abf^+ - \Abf^- |)$. 
This leads us to \eqref{2.33} in the strongly dominant case at least.  

Second, when the $i$-component is dominant but not strongly dominant, we can 
assume (the other case being similar) that  
$$
\kappa_2 \,  | \lam_{i+}^\Abf - \lam_{i-}^\Abf | \, |\alpha_i^+| 
\leq | r_{i+}^\Abf - r_{i-}^\Abf | \, |\beta_i^-| + | \Abf^+ - \Abf^- | \, \sum_{k\neq i} |\beta_k^-|,
$$
thus by \eqref{eigenv} and \eqref{beta} 
\be
{\kappa_2 \over 2} \,  | \lam_{i+}^\Abf - \lam_{i-}^\Abf | \,| \alpha_i^+| 
\leq | \Abf^+ - \Abf^- | \, \sum_{k\neq i} |\beta_k^-|,
\label{2.34} 
\ee
Here, the weight satisfies solely \eqref{contraintes-i2} and we write 
$$ 
\aligned 
\beta_i^- \, \wbf_i^- + \beta_i^+ \, \wbf_i^+
& = \sgn(\beta_i^-) \, ( |\beta_i^-|  - |\beta_i^+| ) \, \wbf_i^-  + \sgn(\beta_i^-) \, |\beta_i^+| \, ( \wbf_i^- - \wbf_i^+ ) 
\\
& \leq  K \, O(1) \, | \lam_{i+}^\Abf - \lam_{i-}^\Abf | \, |\beta_i^-| + O(| \Abf^+ - \Abf^- |) \, \sum_{k\neq i} |\beta_k^-|. 
\endaligned 
$$  
But \eqref{2.34} implies 
$$
\aligned 
K \, O(1) \, | \lam_{i+}^\Abf - \lam_{i-}^\Abf | \, |\beta_i^-|
& \leq {K \over 2} \, | \Abf^+ - \Abf^- | \, \sum_k |\beta_k^-|. 
\endaligned 
$$
provided $\delta_0, \kappa_2$ are sufficiently small.  
In turn, taking into account the other characteristic components which yield a decay 
of $K \, | \Abf^+ - \Abf^- | \, \sum_k |\beta_k^-|$,  we arrive at 
$$  
\aligned 
B 
& \leq 
-K \, | \Abf^+ - \Abf^- | \sum_{k} |\beta_k^-| 
+ {K \over 2} \, | \Abf^+ - \Abf^- | \, \sum_k |\beta_k^-| 
\\ 
&= - {K \over 2} \, | \Abf^+ - \Abf^- | \, \sum_{j \neq i} |\beta_j^-|,
\endaligned 
$$ 
and we conclude that \eqref{2.33} holds for all undercompressive shocks. 

\vskip.3cm 

\noindent{\bf Case 3. }  Finally, consider the case of an $i$-rarefaction-shock, for 
which no constraint has been imposed on the component $\wbf_i$. Here, we will show that 
\be
\aligned 
B \leq 
& C \, | \Abf^+ - \Abf^- | \, |\psi_+ - \psi_-| - {K \over 3}  | \lam_{i+}^\Abf - \lam_{i-}^\Abf | \, |\beta_i^-| 
\\
& - {K \over 3} | \Abf^+ - \Abf^- | \, \sum_{k\neq i} |\beta_k^-| \qquad \text{ in } \Rcal_i^\Abf, 
\endaligned 
\label{2.35} 
\ee
as well as the sharper estimate on non-dominant components: 
\be
\aligned 
B \leq 
& - {K \over 4}  | \lam_{i+}^\Abf - \lam_{i-}^\Abf | \, |\beta_i^-|
   - {K \over 4} | \Abf^+ - \Abf^- | \, \sum_{k\neq i} |\beta_k^-| \qquad \text{ in }  \Rcal_i^\Abf \text{ and non-dominant.}
\endaligned 
\label{2.35b} 
\ee

We first derive \eqref{2.35} for all rarefaction-shocks. From \eqref{2.31} we get 
\be
\aligned 
B & \leq
  2 \, w^{\max} \, |\beta_i^-| + O(| \lam_{i+}^\Abf - \lam_{i-}^\Abf |) \, |\beta_i^-|
  - {K \over 2} \, | \Abf^+ - \Abf^- | \, \sum_{j \neq i} |\beta_j^-| 
  \\ 
& \leq
  3 \, w^{\max} \, |\beta_i^-| - {K \over 2} \, | \Abf^+ - \Abf^- | \, \sum_{j \neq i} |\beta_j^-|  
\endaligned 
\label{2.36}
\ee
and we distinguish between two subcases: 
\vskip.3cm 

-- If $\alpha_i^- \, \alpha_i^+ \geq 0$, then by \eqref{2.13} we have 
$$
(\lamb^\Abf - \lam_{i+}^\Abf) \, \alpha_i^+ + (\lam_{i-}^\Abf - \lamb^\Abf) \, \alpha_i^- 
= O(| \lam_{i+}^\Abf - \lam_{i-}^\Abf |) |\beta_i^-|   + O(| \Abf^+ - \Abf^- |) \sum_{k\neq i} |\beta_k^-|.
$$
The two terms on the left-hand side above have the {\sl same} sign, therefore
$$ 
\aligned 
|\beta_i^+| + | \beta_i^-| 
& = \bigl|(\lamb^\Abf - \lam_{i+}^\Abf) \, \alpha_i^+ \bigr| + \bigl|(\lam_{i-}^\Abf - \lamb^\Abf) \, \alpha_i^- \bigr|  
\\
& = O(| \lam_{i+}^\Abf - \lam_{i-}^\Abf |) |\beta_i^-|   + O(| \Abf^+ - \Abf^- |) \sum_{k\neq i} |\beta_k^-|. 
\endaligned 
$$
We can then suppress the term $O(| \Abf^+ - \Abf^- |)$ above by taking $K$ sufficiently large, 
and so 
$$ 
B 
\leq -{ K \over 3} \, | \Abf^+ - \Abf^- | \, \sum_{j \neq i} |\beta_j^-|,  
$$
which --using once more the previous inequality-- implies \eqref{2.35}.
\vskip.3cm

-- If $\alpha_i^- \, \alpha_i^+ < 0$, then using $\lam_{i+}^\Abf \leq \lamb^\Abf \leq \lam_{i-}^\Abf$ we find 
\be
\aligned 
|\beta_i^-| = |\lam_{i-}^\Abf - \lamb^\Abf| \, |\alpha_i^-| 
& \leq |\lam_{i+}^\Abf - \lam_{i-}^\Abf| \, |\alpha_i^+ - \alpha_i^-| 
\\ 
& \leq O(| \Abf^+ - \Abf^- |) \, |\psi_+ - \psi_-|, 
\endaligned 
\label{2.37}
\ee
so that 
$$ 
\aligned 
B 
&  \leq O(| \Abf^+ - \Abf^- |) \, |\psi_+ - \psi_-|  
- {K \over 2} \, | \Abf^+ - \Abf^- | \, \sum_{j \neq i} |\beta_j^-|, 
\endaligned 
$$
where, in the latter, the estimate \eqref{2.37} on $|\beta_i^-|$ was used 
once more. Since the total variation $\sum_{x \in \Jcal(t)} | \Abf^+(t,x) - \Abf^-(t,x) |$ is uniformly bounded in $t$, 
this proves \eqref{2.35}. 

Finally, to treat the case of a rarefaction-shock that is not strongly dominant, we modify 
the argument in \eqref{2.37} by using the condition \eqref{2.34}:
\be
\aligned 
|\beta_i^-| = |\lam_{i-}^\Abf - \lamb^\Abf| \, |\alpha_i^-| 
& \leq |\lam_{i+}^\Abf - \lam_{i-}^\Abf | \, |\alpha_i^-| 
\\ 
& \leq O(| \Abf^+ - \Abf^- |) \, \sum_{k \neq i} |\beta_k^-|, 
\endaligned 
\label{2.37b}
\ee
which leads to \eqref{2.35b}.

The proof of Proposition~\ref{LI-main} is completed. 
\end{proof}


\section{A new functional for nonlinear hyperbolic systems}
\label{WE-0} 

\subsection{Scalar conservation laws}
 
Conditions on the weight function $\wbf$ were proposed in Section~\ref{section4}, 
and our aim is now to show that such a weight exists. 
More precisely, we introduce now slightly different conditions
which are easier to work with at this stage of the analysis
and will be shown in the following section to imply the conditions required for our analysis of Section~\ref{section4}.  
Importantly, we are going to see now 
that wave cancellations must be taken into account to establish that the weight remains uniformly bounded. 
The new technique proposed here relies on the adjoint system \eqref{2.adjoint} 
and on the property of generalized characteristics for the averaged matrix.  
The following is a natural generalization of the technique proposed earlier 
for genuinely nonlinear equations \cite{HuLeFloch,LeFloch-book} which also rely on wave cancellation. 

For scalar equations, we use Dafermos' front tracking technique \cite{Dafermos-front}.  

\begin{theorem}[Functional for scalar conservation laws]
\label{51} 
Consider a scalar conservation law with general flux $f$. 
Given $C_1, K>0$ there exists $C_2=C_2(C_1,K)$ such that the following property holds. 
Given (piecewise constant) front tracking solutions $u_h,u_h'$ satisfying the uniform bounds 
\be
\label{Cone}
\| u_h \|_{L^\infty} + \| u_h' \|_{L^\infty} + \sup_t TV(u_h(t)) + \sup_t TV(u_h'(t)) \leq C_1, 
\ee
there exists a (piecewise constant) weight function $\wbf = \wbf[u_h,u_h']$ satisfying 
$$
1/C_2 \leq \wbf \leq C_2
$$  
together with the following conditions on undercompressive shocks of the averaged speed $\abf:=\ab(u_h,u_h')$: 
\be
\label{contr-scalar}
\wbf^+ - \wbf^-   = 
\begin{cases}  
- K \, | \abf^+ - \abf^- |,  & \Scal^\abf, 
\\ 
\hskip.3cm  K \, | \abf^+ - \abf^- |,  & \Fcal^\abf. 
\end{cases}
\ee 
\end{theorem}

\

From this theorem, we deduce that the corresponding weighted $L^1$ functional is (essentially) decreasing: 
\be
\label{dissipa} 
\aligned 
& \int_\R |u_h' - u_h | \, \wbf[u_h,u_h'](t_2) \, dx 
+ \int_{t_1}^{t_2} \sum_{\Lcal^\abf}  |\lamb^\abf - \abf^- | \, | u_- - u_-'| \, dt
\\
& + K \, \int_{t_1}^{t_2} \sum_{\Scal^\abf \cup \Fcal^\abf}  |\abf^+ - \abf^- | \, |\lamb^\abf - \abf^- | \, | u_- - u_-'| \, dt 
\\
&  \leq \int_\R | u_h' - u_h | \, \wbf[u_h,u_h'](t_1) \, dx + O(h), \qquad t_2 \geq t_1, 
\endaligned 
\ee 
where $\lamb^\abf$ denotes the shock speed of the function $\abf$. 
The remainder $O(h)$ is bounded by the maximum size of rarefaction fronts in $u_h,u_h'$, which 
tends to zero with the discretization parameter $h$. Hence, by letting $h \to 0$ in \eqref{dissipa} one deduces  
that the limits $u:= \lim u_h$ and  $u':= \lim u_h'$ satisfy the sharp $L^1$ contraction property 
$$
\int_\R |u' - u| \, \wbf(t_2) \, dx < \int_\R |u' - u| \, \wbf(t_1) \, dx,  
\qquad t_2 > t_1,
$$
where $\wbf := \lim \wbf[u_h,u_h']$ 
and the inequality is strict (except in the trivial case $u \equiv u'$). Moreover, sharp dissipation terms can 
be made explicit, as was done for convex flux 
by Dafermos \cite{Dafermos-book} and Goatin and LeFloch \cite{GoatinLeFloch-scalar}.

\begin{proof} 
1. For simplicity in the discussion and without loss of generality (modulo arbitrarily small perturbations 
of the data), we can assume that all shocks in the function $\abf$ 
satisfy the inequalities in Definition~\ref{def21} in a {\sl strict} sense and that 
$ \psi := u_h - u_h'$ does not vanish, but of course may change sign at jump points.
(This is possible since $u_h,u_h'$ are piecewise constant with finitely many waves only.) 

Decompose the space-time $\R_+\times\R$ into finitely many maximal regions $\Omega$ where $u_h - u_h'$ keep a constant sign
and which, therefore, are limited by polygonal curves of changes of sign for the function $\psi$. 
According to Lemma~\ref{Lemma3.2} the boundaries of $\Omega$ consist of compressive or rarefaction shocks
along which no constraint is imposed in \eqref{contr-scalar}. Consequently, we can turn attention to defining 
the weight $\wbf = \wbf(t,x)$ within a given region $\Omega$.  

2. Consider backward generalized characteristics $y_m : [\bart_m, \tbar_m] \subset [0,+\infty) \to \R$ 
associated with the averaged speed $\abf$ and originating and finishing on the boundary of $\Omega$, 
$$
{d y_m\over dt} (t) \in I(\abf_-(t,y_m(t)), \abf_+(t,y_m(t)) \big),  \quad t \in [\bart_m, \tbar_m], 
$$
where $I(b,c) := (\min(b,c), \max(b,c))$. 
Observe that backward generalized characteristics {\sl within} $\Omega$ are unique once their origin is given, since 
we assumed that no shock in $\abf$ is left- nor right-characteristic. 
 
More precisely, since $u,u'$ admit finitely many waves and interaction points, we can select 
a {\sl complete} family of generalized characteristics $y_m$, $m=1, \ldots, M$ within the region $\Omega$, in the 
sense that every other characteristic meets precisely the same waves of $u_h,u_h'$ as one of the curves $y_m$. 
 
\
 
3. We now determine the weight within the region $\Omega$ by {\sl formally} solving the adjoint equation  
\be
\label{eq} 
\del_t \wbf +  \abf \, \del_x \wbf = 0 \quad \text{ away from shocks,} 
\ee 
in the following way. We prescribe the weight at the final point $\tbar_m$ of each characteristic
(when $\tbar_m<\infty$) 
$$
\wbf (\tbar_m,y_m(\tbar_m)) = C_0, \quad m=1, \ldots, M,
$$
where $C_0$ is a (large) positive constant. When $\tbar_m<\infty$ we require that $\wbf (t,y_m(t))$ equals 
$C_0$ for all sufficiently large times $t$. 
Then, we determine the value $\wbf(t,y_m(t))$ of the weight at arbitrary times $t \in (\bart_m, \tbar_m)$
by following the generalized characteristic $y_m$ backward in time from the final time $\tbar_m$
and by requiring that the weight:
\par 
--  remains constant away from shocks (so that \eqref{eq} holds), and 
\par
-- jumps up or down, according to the constraint \eqref{contr-scalar}, 
when it crosses a shock.   
\par
\noindent  

Hence, our construction determines a weight defined along the curves $y_m$ and, in summary, we can write  
(away from shocks) 
\be
\label{sum}
\wbf(t,y_m(t)) = C_0 + K \sum_{s \in [t,\tbar_m)} \eta_m(s) \, \big( \abf^+(s,y_m(s)) - \abf^-(s,y_m(s)) \big),
\ee
where the sign $\eta_m(s) = \pm$ is determined from \eqref{contr-scalar} and depend whether the shock
is slow or fast undercompressive and whether it is crossed from left to right or from right to left. 

In turn, since the curves $y_m$ is a complete family of characteristics, the weight can be uniquely extended
as a piecewise constant function $\wbf = \wbf(t,x)$ defined in the whole region $\Omega$. 

\

4. Next, for any fixed time $T>0$ we can introduce wave partitions (Liu \cite{Liu2}) of the waves in 
both solutions $u_h,u_h'$ within the slab $[0,T] \times \R$. 
By definition, the solution $u_h$ can be regarded as the superposition of finitely many 
{\sl elementary waves} $(u_p, u_{p+1})$ associated with piecewise affine functions $\varphi_p : [0,T_p] \to \R$. 
The latter determine the trajectories of the elementary waves in the plane. The wave $(u_p, u_{p+1})$
is born at the initial time $t=0$ and completely cancelled out at the time $t=T_p$. The cancellation takes place at an interaction 
point with another wave. In the scalar case, no new wave is generated at interactions, and waves can only be split or cancelled out. 
A similar notation, $(u_q', u_{q+1}')$ and $\varphi_q' : [0,T_q'] \to \R$ is used for the solution $u_h'$.

5. We now establish that the function $\wbf$ remains uniformly bounded. It is (necessary and) 
sufficient to bound the maximum oscillation of $\wbf$ within the region $\Omega \subset [0,T] \times \R$
for arbitrary $T$ and along each characteristic $y_m$. We will prove that 
\be
\label{oscill}
\aligned 
\Osc_m 
& := \sup_{t \in  (\bart_m, \tbar_m)} \wbf(t,y_m(t)) - \inf_{t \in  (\bart_m, \tbar_m)} \wbf(t,y_m(t))
\\
& \lesssim TV(u_h(0)) + TV(u_h'(0)),
\endaligned 
\ee
where the implied constant depends only upon $C_1,K$ arising in \eqref{Cone}-\eqref{contr-scalar}. 

Recall first that $\Omega$ contains only undercompressive shocks and that, within the region $\Omega$, the 
generalized characteristic $y_m$ crosses undercompressive shocks {\sl transversally.} 
To derive \eqref{oscill} we need to take into account {\sl wave cancellations} and to rely on the following
{\sl linearity property} of the averaged speed $\ab$ with respect to elementary waves: 
if a wave $(u_-,u_+)$ in the solution $u_h$ is split into elementary waves $(u_p, u_{p+1})$ then for every constant $u'$:
\be
\label{linear}
\ab(u_+,u') - \ab(u_-,u') = \sum_p \ab(u_{p+1}, u') - \ab(u_p,u'), 
\ee
which we refer to as the linearity property. (Of course, this property holds for arbitrary functions.)

To estimate the weight $\wbf$ along a path $y_m$ we fix an interval $[t_1,t_2] \subset (\bart_m, \tbar_m)$
and consider the set of all elementary waves $(u_p, u_{p+1})$, $p \in E=E_m(t_1,t_2)$, 
and $(u_q', u_{q+1}')$, $q \in F=F_m(t_1,t_2)$, that are crossed by the characteristic $y_m$ at least once
(but possibly more than once). Observe that a slow shock can only be crossed from right to left as we move backward, 
while the opposite is true for a fast shock.  At this juncture we need to recall Lemma~\ref{Lemma3.2} which implies that
a specific sign $\eta_p = \pm$ can be associated with a given elementary  wave $(u_p,u_{p+1})$
which only depend whether $u_h - u_h'$ is positive or negative in the region $\Omega$, 
whether the wave is increasing or decreasing. (Note in passing that signs are different also for the waves $(u_p,u_{p+1})$
and the waves $(u_p',u_{p+1}')$.) Then, at a point of discontinuity $(s,y_m(s))$
 the sign $\eta_m(s)$ arising in \eqref{sum} is $\pm \eta_p$ 
and {\sl alternatively positive and negative} as the wave is crossed from left to right and then 
move back to the left-hand of the characteristic at a later time. 

Consequently, given a pair $(p,q) \in E \times F$, together with the trajectories $\varphi_p, \varphi_q'$
the contribution of the waves  $(u_p,u_{p+1})$ and  $(u_q',u_{q+1}')$ to the weight along the curve $y_m$
can be estimated by the strengths of these waves
$$
|u_{p+1} - u_p| + | u_{q+1}' - u_q' |
$$
plus the alternating sum 
$$
\aligned 
& \Sigma_{m,pq}(t_1,t_2)
\\
&  := \sum \pm \Big( \big(\ab(u_{p+1}, u_{q+1}') - \ab(u_p, u_{q+1}')\big)
 - \big( \ab(u_{p+1}, u_q') - \ab(u_p, u_q') \big) \Big),  
\endaligned 
$$
over all pairs of waves crossings with the characteristic $y_m$ within the time interval $[t_1,t_2]$. 
The sign is alternately positive and negative as one moves along the characteristic. 
Consequently, we can keep only one term at most 
$$
\aligned 
|\Sigma_{m,pq}(t_1,t_2)| 
& \leq  \big|\ab(u_{p+1}, u_{q+1}) - \ab(u_p, u_{q+1}') - \ab(u_{p+1}, u_{q+1}) - \ab(u_p, u_{q+1}') \big|
\\
& \lesssim | u_{p+1} - u_p | \, | u_{q+1}' - u_q'|. 
\endaligned 
$$

In turn, by taking into account all pairs $(p,q)$ of elementary waves we conclude that 
$$
\aligned 
\Osc_m 
& \lesssim \sum_p | u_{p+1} - u_p |  + \sum_q | u_{q+1}' - u_q'|
    + \sum_{p,q} | u_{p+1} - u_p | \, | u_{q+1}' - u_q'|
\\
& \leq TV(u_h(0)) + TV(u_h'(0)) + TV(u_h(0)) \, TV(u_h'(0)),  
\endaligned 
$$
which under the assumptions \eqref{Cone} implies \eqref{oscill}.
\end{proof}


\subsection{Nonlinear hyperbolic systems} 
\label{522}

We now turn attention to general systems \eqref{1.1} 
and to piecewise constant solutions $u_h$ generated by wave front tracking. 
For the actual construction and properties of front tracking approximation\footnote{See the recent preprint: O. Glass and P.G. LeFloch, 
Nonlinear hyperbolic systems:
Non-degenerate flux, inner speed variation, and graph solutions, Arch. Rational Mech. Anal., to appear.}
we refer to the general theory in \cite{Bianchini,IguchiLeFloch,LiuYang-exist} 
(for general systems via the Glimm scheme) as well as to the 
earlier references \cite{Dafermos-front,DiPerna,Risebro,Bressan1,Bressan2,LeFloch-book}  
(for genuinely nonlinear systems via front-tracking). 

Specifically, given two approximate solutions $u_h,u_h'$ we introduce the averaged matrix $\Abf = \Ab(u_h,u_h')$
and the function $\psi = u_h' - u_h$. Then, by \eqref{charact} we define the characteristic components $\alpha_j$
and their traces $\alpha_j^\pm=\alpha_j(t,x\pm)$ at a shock $(t,x)$. Without loss of generality we can assume that shock trajectories in $u_h, u_h'$ cross at finitely many points only and never superimpose on an open interval of time, 
and that the components $\alpha_j$ never vanish while the shocks are non-characteristic on both sides. 
This can be achieved by arbitrarily small perturbations of the data. Following \cite{HuLeFloch}
we define the strength $\eps^\Abf = \eps^\Abf(t,x)$ at a shock of the matrix-valued map $\Abf$ as 
the strength of the corresponding shock in the solutions $u_h, u_h'$. The latter  
is measured with the global parameter introduced in Section~\ref{subsection3.3}: 
$$
\eps^\Abf =\eps^\Abf(t,x) := \begin{cases}
|\mu_i(u_h^+) - \mu_i(u_h^-)|, \quad \text{ $i$-shock in } u_h, 
\\
|\mu_i({u_h'}^+) - \mu_i({u_h'}^-)|,  \quad \text{ $i$-shock in } u_h', 
\end{cases}
$$
Clearly, we have 
$$
|\Abf^+ - \Abf ^-| \lesssim \eps^\Abf \quad \text{ at shocks.}   
$$ 

We now introduce conditions on the weight $\wbf=\wbf(t,x)$ which are most convenient 
to work with at this stage. Recall first that the relevance of the signs of the characteristic components
has been emphasized in \eqref{signes-j}-\eqref{signes-i}, and in Lemma~\ref{Prop3.5}. 
In particular, the latter allows us to decompose the set of all $i$-discontinuities in $u_h, u_h'$ into two sets 
according to whether the quantity (defined at a jump point $(t,x)$ from the left- and right-hand traces) 
$$
\rho_i(t,x) := \rho_i(u_h^-,u_h^+,{u_h'}^-,{u_h'}^+)
=
\begin{cases}
\rho_i(u_h^-,u_h^+,u_h'), \qquad \text{ shock in } u_h, 
\\
\rho_i({u_h'}^-, {u_h'}^+,u_h),  \quad \text{ shock in } u_h', 
\end{cases}
$$
is negative or positive respectively.  
$$
\Jcal_i^\Abf 
                  = \big\{ \rho_i(t,x) <0 \big\} \cup  \big\{ \rho_i(t,x) >0 \big\}. 
$$
Lemma~\ref{Prop3.5} shows some compatibility between 
the signs of the jump of the weight $\wbf_i^+ - \wbf_i^-$ 
and the jump of the eigenvalue 
$\lamb_i(u_h^+, {u_h'}^+) - \lamb_i(u_h^-, {u_h'}^-) =  \lam_{i+}^\Abf  - \lam_{i-}^\Abf$. 

Actually, the above is true for strongly dominant shocks only. We propose here to relax this condition 
and to impose the conditions suggested by Lemma~\ref{Prop3.5} {\sl even for} 
non-strongly dominant 
components. We will see that this strategy leads to a well-defined and well-behaved weight $\wbf$. 

In other words, we require that  for all $j \neq i$ 
\be
\label{cond1}
\wbf_j^+ - \wbf_j^- = 
\begin{cases}  
\hskip.3cm K \, \eps^\Abf,  & j < i   \text{ and } \alpha_j^- \, \alpha_j^+ > 0,  
\\ 
- K \, \eps^\Abf,                   & j > i \text{ and }  \alpha_j^- \, \alpha_j^+ > 0, 
\end{cases}
\ee   
while for $j=i$ and all undercompressive $i$-shocks
\be
\label{cond2}
\wbf_i^+ - \wbf_i^- = 
\begin{cases}  
- K \, ( \lam_{i+}^\Abf  - \lam_{i-}^\Abf ),  &  \rho_i <0  \text{ and } \alpha_i^- \, \alpha_i^+ > 0, 
\\ 
\hskip.3cm  K \,  ( \lam_{i+}^\Abf  - \lam_{i-}^\Abf ),   & \rho_i >0  \text{ and } \alpha_i^- \, \alpha_i^+ > 0,
\end{cases}
\ee 
which corresponds to well-defined conditions for each increasing/decrea\-sing shocks in regions where 
$\alpha_i$ is positive/negative. 
The proof of the following theorem will rely on a wave partition into increasing/decreasing elementary waves.

\begin{theorem}[$L^1$ functional for systems of conservation laws]
\label{weight-sys} 
Consider a strictly hyperbolic system of conservation laws with general flux $f$. 
Given $C_1, K>0$ there exists $C_2=C_2(C_1,K)$ such that the following property holds. 
Given (piecewise constant) front tracking solutions $u_h, u_h'$ satisfying the uniform bounds 
$$
\| u_h \|_{L^\infty} + \| u_h' \|_{L^\infty} + \sup_t TV(u_h(t)) + \sup_t TV(u_h'(t)) \leq C_1, 
$$
there exists a (piecewise constant) weight function $\wbf = (\wbf_j[u_h,u_h'])_{1 \leq j \leq N}$ satisfying 
$$
1/C_2 \leq \wbf_j[u_h,u_h']  \leq C_2
$$ 
together with the conditions \eqref{cond1} and \eqref{cond2}.  
\end{theorem}

\begin{proof} 
1. Our construction will be similar to the one proposed for scalar equations, with however
some important modifications due to 
\par 
-- the new waves generated at interactions, and 
\par
-- the nature of shocks (fast/slow, compressive/rarefaction) which is not entirely determined by the signs of 
the characteristic components $\alpha_j^\pm$ and the coefficient $\rho_i$. Indeed, 
Lemma~\ref{Prop3.5} is concerned only with dominant components $\alpha_j$ ($j \neq i$)
and strongly dominant components $\alpha_i$ across $i$-shocks.   
\par 
\noindent However, by taking into account contributions due to small waves 
the construction can still be carried out and wave cancellations be taken into account, as we now explain. 
Fixing some index $j$ we focus on defining the component $\wbf_j$ of the weight. 
 
\

2. Decompose $\R_+\times\R$ into maximal regions $\Omega_j$ 
where the characteristic component $\alpha_j$ keeps a constant sign
and which are limited by polygonal curves across which $\alpha_j$ changes sign. 
According to Lemma~\ref{Prop3.5} the boundary of a region $\Omega_j$ consists of 
\par 
-- $k$-shocks ($k\neq j$) with dominant characteristic components $\alpha_j^\pm$ (but opposite signs), or
\par 
-- compressive or rarefaction $j$-shocks with dominant characteristic components $\alpha_j^\pm$
(with opposite signs), or 
\par 
-- shocks with non-dominant components $\alpha_j^\pm$ (with opposite signs), 
\par 
-- or possibly an interval of the initial line $t=0$. 
\par 
\noindent 
We observe that along of the first type of boundary the simpler constraint  \eqref{cond1} is imposed while 
along the other types of boundary no constraint is imposed on the $j$-component of the weight $\wbf_j$. 

Consequently, we need not only to deal with $\wbf_j$ in one of the regions $\Omega_j$
(in which the constraints \eqref{cond1}-\eqref{cond2} are relevant),
but also to take into account waves leaving one of the regions $\Omega_j$.

Moreover, in contrast with the case of scalar equations, a region $\Omega_j$ may now contain not only 
undercompressive $j$-shocks but also 
\par
-- $k$-shocks ($k\neq j$) and 
\par 
-- compressive shocks or rarefaction-shocks of the $j$-family (which however must be non-dominant). 

\

3. We will rely on a sufficiently large family backward generalized $j$-charac\-teristics $y_m^j : [0,\infty) \to \R$ 
associated with the eigenvalue $\lam_j^\Abf$:   
$$
{d y_m^j \over dt} (t) \in I\big( \lam_{j-}^\Abf(t,y^j_m(t)), \lam_{j+}^\Abf(t,y^j_m(t)) \big),  \quad t \in (0,\infty). 
$$
Here, a given region $\Omega_j$ may contain compressive $j$-shocks
so that backward generalized characteristics need not be unique.  
Also, such regions may well allow $j$-waves to cross their boundary: we impose that $j$-characteristics
never exist a region $\Omega_j$, except along a boundary consisting of a $k$-wave with $k \neq j$. 
Therefore $y_m^j$ may only exit to the right of $\Omega_j$ by crossing a $k$-wave with $k < j$, 
or exit to the left by crossing a $k$-wave with $k>j$. And a $j$-characteristic does not cross  
compressive and rarefaction $j$-shocks, but may travel with the same location and speed as a rarefaction shock
on some time interval. 

Since the solutions $u_h,u_h'$ admit finitely many waves and interaction points, we can select a {\sl complete} family 
of generalized $j$-characteristics $y_m$, $m=1, \ldots, M$ covering the whole of $\R_+ \times \R$,
such that every other $j$-characteristic meets precisely the same waves of $u,u'$ as one of the curves 
in the family $y_m^j$. 

\
 
4. Next, backward from the time $T>0$ we can determine wave partitions within the slab $[0,T] \times \R$
of the waves in both solutions $u_h,u_h'$ (Liu \cite{Liu2}). 
By definition of a wave partition, the solution $u_h$ is regarded as the linear superposition of finitely many 
{\sl elementary $k$-waves} $(u_p^k, u_{p+1}^k)$ with trajectories $\varphi_p^k : [\barT_p^k,\Tbar_p^k] \to \R$. 
The wave $(u_p^k, u_{p+1}^k)$ is born at the initial time $\barT_p^k$ and is completely cancelled at the time $t=\Tbar_p^k$. 
For systems, new waves may be generated at wave interaction points.  
The total strength, and the change along their trajectories, of completely cancelled waves ($\Tbar_p^k < T$), 
new waves  ($0 < \barT_p^k $),  and surviving waves $\Tbar_p^k = T$)  
can be estimated uniformly by the initial total variation, the interaction potential, and the cancellation measure, 
all of them being in turn controled by the total variation of the initial data $u_h(0)$.  
A similar notation, $({u'}_q^k, {u'}_{q+1}^k)$ and ${\varphi'}_q^k : [0,T_l'] \to \R$ will be used for the 
elementary $k$-waves in the solution $u_h'$. 

\ 

5. We now determine the weight by {\sl formally} solving the adjoint equation  
\be
\label{eq22} 
\del_t \wbf +  \Abf \, \del_x \wbf = 0 \quad \text{ away from shocks,} 
\ee 
in the following way. We prescribe the weight at some sufficiently large time $T$ along each characteristic
$$
\wbf_j (T,y_m^j(T)) = C_0, \quad m=1, \ldots, M,
$$
where $C_0$ is a (large) positive constant.  
Then, we determine the value $\wbf_j(t,y_m^j(t))$ of the weight at arbitrary times $t < T$
by following the $j$-characteristic $y_m^j$ backward in time and requiring that $\wbf_j$ 
\par 
--  remains constant away from shocks (so that \eqref{eq22} holds), 
\par
-- jumps up or down, according to the constraint \eqref{cond1} when it crosses a $k$-shock ($k \neq j$), and 
\par 
-- jumps up or down, according to the constraint \eqref{cond2} when it crosses an undercompressive $j$-shock. 
\par
\noindent By construction compressive and rarefaction $j$-shocks are never crossed by $j$-cha\-racteristics. 

Our construction determines a function defined along each of 
the curves $y_m^j$, and away from shocks we can write  
\be
\label{sum22}
\aligned 
\wbf_j(t,y_m^j(t)) = C_0 
& + K \sum_{s \in [t,T) \atop \text{$j$-shocks}} \eta_m(s) \, \big( \lam_{j+}^\Abf(s,y_m^j(s)) - \lam_{j-}^\Abf(s,y^j_m(s)) \big),
\\
& + K \sum_{s \in [t,T) \atop \text{$k$-shocks}}  \eta_m(s) \, \eps^\Abf(s,y_m^j(s)),  
\endaligned 
\ee
where the sign $\eta_m(s) = \pm$ is determined from \eqref{cond1}-\eqref{cond2} and depends whether the shock is 
slow or fast undercompressive and whether it is crossed from left to right or from right to left. 

In turn, by using that the curves $y_m^j$ form a complete family of characteristics, the weight is uniquely extendable 
as a piecewise constant function $\wbf_j = \wbf_j(t,x)$ defined in $\R_+\times\R$. 

\

6. We now claim that the function $\wbf_j$ remains uniformly bounded, independently of the parameter $h$. 
We need to estimate the maximum oscillation of $\wbf_j$  
in $[0,T] \times \R$ along each $j$-characteristic $y^j_m$. We will actually prove that 
\be
\label{oscill22}
\aligned 
\Osc_m^j  
& := \sup_{t \in [0,T]} \wbf_j(t,y^j_m(t)) - \inf_{t \in [0,T]} \wbf_j(t,y^j_m(t))
\\
& \lesssim TV(u_h(0)) + TV(u_h'(0)). 
\endaligned 
\ee 
We follow here the argument from the scalar case, but now the regions $\Omega_j$ 
may contain not only undercompressive shocks 
(which the generalized characteristic $y_m$ crosses transversally) but also $k$-shocks and other (non-dominant) $j$-waves. 
 
The {\sl linearity property} with respect to elementary waves, pointed out for scalar equations, is 
now expressed in terms of the $j$-eigenfunction $\lamb_j = \lam_j^\Abf$: if a $j$-wave $(u_-,u_+)$ in the solution $u_h$ is split into 
elementary waves $(u_p^j, u_{p+1}^j)$ then for every constant $u'$:
\be
\label{linear666}
\lamb_j(u_+,u') - \lamb_j(u_-,u') = \sum_p \lamb_j(u_{p+1}^j, u') - \lamb_j(u_p^j,u'). 
\ee 

To estimate the weight $\wbf_j$ along a given path $y_m^j$ we consider an arbitrary interval $[t_1,t_2] \subset [0,T]$
and consider the set of all elementary waves crossed by the characteristic $y_m^j$ at least once
(but possibly more than once), precisely: 
\par 
-- all elementary $j$-waves $(u_p^j, u_{p+1}^j)$, $p \in E^j_m = E_m^j(t_1,t_2)$, associated with the solution $u_h$, 
\par 
-- all elementary $k$-waves $(u_p^k, u_{p+1}^k)$, $p \in E^k_m = E_m^k(t_1,t_2)$,  with $k \neq j$, 
\par 
-- all elementary $j$-waves $({u_q'}^j, {u'}_{q+1}^j)$, $q \in F_m^j=F_m(t_1,t_2)$, associated with the solution $u_h'$, and
\par 
-- all elementary $k$-waves $({u_q'}^k, {u'}_{q+1}^k)$, $q \in F_m^j=F_m(t_1,t_2)$, associated with the solution $u_h'$.

Observe that a slow $j$-shock can only be crossed from right to left as we move backward along a $j$-characteristic, 
while we have the opposite direction for a fast shock. 
Relying on the key observation in Lemma~\ref{Prop3.5} we observe that a definite sign 
$\eta_p^l = \pm$ (depending whether the wave is increasing or decreasing)
can be associated with each elementary  wave $(u_p^j,u_{p+1}^j)$, so that 
the sign $\eta_m(s)$ in \eqref{sum22} is $\pm \eta_p^l$ where the sign $\pm$ only depends 
\par 
-- whether the characteristic $y_m^j$ at the time $s$ lies in a region $\Omega_j$ where the characteristic
component is positive and $-\eta_p^l $ in a region where it is negative, 
\par 
-- whether the wave is crossed from left to right, or from right to left. 
 
Consequently, as the characteristic $y_m^j$ crosses through the elementary waves of $u_h, u_h'$, 
{\sl alternatively positive and negative} signs arise in the first sum in \eqref{sum22}.
Indeed, this is true in a given region $\Omega_j$, as well as when the wave crosses its 
boundary and enter another region of constant sign.

Consequently, given a pair $(p,q) \in E_m^j \times F_m^j$, together with the trajectories $\varphi_p^j, {\varphi'}_q^j$
the contribution of the $j$-waves  $(u_p^j,u_{p+1}^j)$ and  $({u'}_q^j,{u'}_{q+1}^j)$ to the weight $\wbf_j$ 
along the curve $y_m^j$ within the time interval $(t_1,t_2)$ is estimated by 
$$
|u_{p+1}^j - u_p^j| + | {u'}_{q+1}^j - {u'}_q^j |
$$
plus the sum 
$$
\aligned 
& \Sigma^j_{m, pq}(t_1,t_2)
\\
&  := \sum \pm \Big( \big( \lamb_j(u_{p+1}^j, {u'}_{q+1}^j) 
- \lamb_j(u_p^j, {u'}_{q+1}^j)\big) - 
\big(\lamb_j(u_{p+1}^j, {u'}_q^j) - \lamb_j(u_p^j, {u'}_q^j) \big) \Big),  
\endaligned 
$$
over all pairs of crossings of the characteristic $y_m$ within the time interval $[t_1,t_2]$.  
The signs alternate as a given wave passes from left to right, or then back 
from the right to the left of the characteristic, 
and we conclude that terms cancel out two at a time and 
a single term only must be kept at most 
$$
\aligned 
& |\Sigma^j_{m,pq}(t_1,t_2)| 
\\
& \leq 
\Big| \lamb_j(u_{p+1}^j, {u'}_{q+1}^j) - \lamb_j(u_p^j, {u'}_{q+1}^j) - \lamb_j(u_{p+1}^j, {u'}_q^j) + \lamb_j(u_p^j, {u'}_q^j) \Big|,  
\\
& \lesssim |u_{p+1}^j - u_p^j| \, | {u'}_{q+1}^j - {u'}_q^j |. 
\endaligned 
$$

In turn, by taking into account all pairs $(p,q)$ of elementary waves and by making use of 
properties of wave partitions (uniform bounds on cancellation and interaction measures) we find 
$$
\aligned 
& \sum_{p \in E_m^j} |u_{p+1}^j - u_p^j| + \sum_{q \in  F_m^q} | {u'}_{q+1}^j - {u'}_q^j |
\\
& \quad + \sum_{p,q \in E_m^j \times F_m^q} |u_{p+1}^j - u_p^j| \, | {u'}_{q+1}^j - {u'}_q^j |
\\
& \lesssim TV(u_h(0)) + TV(u_h'(0)) + TV(u_h(0)) \, TV(u_h'(0)). 
\endaligned 
$$

On the other hand, to handle the second sum in \eqref{sum22} we observe that the $j$-characteristic $y_m$ 
crosses a $k$-shock with $k \neq j$ at most once, and have suitable alternating properties.
This second sum accounts for the contribution 
$$
\aligned 
& \sum_{p \in E_m^j
\atop k \neq j} |u_{p+1}^k - u_p^j| + \sum_{q \in  F_m^q \atop k \neq j} | {u'}_{q+1}^k - {u'}_q^k |
\\
& \lesssim TV(u_h(0)) + TV(u_h'(0)). 
\endaligned 
$$
The time interval $(t_1,t_2)$ being arbitrary we conclude that 
$$
\aligned 
\Osc_m^j  \lesssim TV(u_h(0)) + TV(u_h'(0)) + TV(u_h(0)) \, TV(u_h'(0)),  
\endaligned 
$$
which shows \eqref{oscill22} and completes the proof of Theorem~\ref{weight-sys}. 
\end{proof}


\section{Continuous dependence property} 
\label{Section6}

We are now in a position to prove the $L^1$ estimate \eqref{1.2}. 
For definiteness we state our results for solutions constructed by front tracking (see Section~\ref{522}). 
With some modification, our method also applies to solutions constructed by the Glimm scheme.

\begin{theorem}[$L^1$ continuous dependence of solutions] 
\label{theorem-main} 
Consider a strictly hyperbolic system \eqref{1.1} with general flux $f$,  
under the assumption that there exists an averaged matrix $\Ab$ satisfying \eqref{matrix0} and 
\be
\label{eigenv2}
| \rb_i(u_+,u') -  \rb_i(u_-,u') | \lesssim | \lamb_i(u_+,u') -  \lamb_i(u_-,u') |, 
\ee
for all states $u_\pm, u'$ under consideration. 
Then, wave front tracking approximations $u_h, u_h'$ with sufficiently small amplitude and total variation
satisfy, for all times $t \geq 0$,  
\be
\|u_h'(t) - u_h(t)\|_{L^1(\R)} \lesssim \|u_h'(0) - u_h(0)\|_{L^1(\R)} + o(h) 
\label{5.4} 
\ee
and, consequently, the limit solutions $u = \lim_{h \to 0} u_h$
and $u' := \lim_{h \to 0} u_h'$ satisfy the $L^1$ continuous dependence property
\be
\|u'(t) - u(t)\|_{L^1(\R)} \lesssim \|u'(0) - u(0)\|_{L^1(\R)}.  
\label{5.5} 
\ee
\end{theorem}
 
 \

This is a generalization of a theorem established earlier  \cite{HuLeFloch, LeFloch-book} 
for the class of genuinely nonlinear systems. Indeed, for such systems, the condition \eqref{eigenv2}
always holds since, on one hand, 
$$
| \rb_i(u_+,u') -  \rb_i(u_-,u') | \lesssim |u_+ - u_-|
$$
and, on the other hand by Lemma~\ref{Lemma3.3},  
$$ 
|u_+ - u_-| \lesssim | \lamb_i(u_+,u') -  \lamb_i(u_-,u') |, 
\quad \text{ genuinely nonlinear fields.} 
$$

We conclude with: 

\begin{corollary} 
\label{corol}
The continuous dependence properties \eqref{5.4} (approximate solutions)
and \eqref{5.5} (exact solutions) hold for  
\begin{enumerate} 
\item the Lagrangian $p$-system of fluid dynamics, and
\item the Euler equations of isentropic fluid dynamics,
\end{enumerate} 
with general pressure-density equations of state. 
\end{corollary}

\begin{proof}[Proof of Theorem~\ref{theorem-main}]
We rely on the framework developed in previous sections. 
Let $\wbf = (\wbf_j(t,x))_{1 \leq j \leq N}$ be the weight determined in Theorem~\ref{weight-sys}. 
We need to check that the conditions \eqref{cond1}-\eqref{cond2} 
imply the conditions \eqref{contraintes-i1}-\eqref{contraintes-i2} required for the stability theory of 
Section~\ref{section4}. 
Once this is checked, the desired $L^1$ estimate \eqref{5.4} follows from 
Proposition~\ref{LI-main} ($L^1$ stability for linear systems) 
and Theorem~\ref{general} (non-existence of strongly dominant rarefaction-shocks), 
and the properties of front tracking approximations, especially the fact that rarefaction-shocks in $u_h, u_h'$ 
have maximal strength converging to zero with $h$. 
 
To compare the conditions in Sections~\ref{section4} and \ref{WE-0}, we 
use the classification in Definition~\ref{def21} and 
the characterization (obtained in Lemma~\ref{Prop3.5}) of the shocks in terms of:
\par
-- the sign of characteristic components $\alpha_j$ and 
\par
-- the sign of the jump of the eigenvalues of $\Abf=\Ab(u_h, u_h')$. 

First of all, in view of \eqref{signes-j} ($\alpha_j^- \alpha_j^+ >0$ for dominant $j$-components of an $i$-shock) 
the conditions \eqref{cond1} imply \eqref{contraintes-j}.  

Second, concerning strongly dominant $i$-components $\alpha_i$, 
we see that \eqref{cond2} in combination with \eqref{3.13} show that 
\be
\label{contraintes-i17}
\wbf_i^+ - \wbf_i^- = 
\begin{cases}  
 - K \,  | \lam_{i+}^\Abf -  \lam_{i-}^\Abf |,  & \Scal_i^\Abf \text{ and strongly dominant,} 
\\ 
\hskip.3cm  K \, | \lam_{i+}^\Abf -  \lam_{i-}^\Abf |, 
&  \Fcal_i^\Abf  \text{ and strongly dominant,} 
\end{cases}
\ee 
which implies \eqref{contraintes-i1}.  

Third, for an undercompressive $i$-shock whose $i$-component is dominant not strongly dominant, 
the sign of the jump $\lam_{i+}^\Abf - \lam_{i-}^\Abf$ is not well defined,  
and the sign of the jump $\wbf_i^+ - \wbf^-$ may not agree with the general rule \eqref{updown}. 
However, from \eqref{cond2} we can still deduce that 
\be
\label{contraintes-i29}
|\wbf_i^+ - \wbf_i^-| =  
K \,  | \lam_{i+}^\Abf - \lam_{i-}^\Abf |, \qquad 
 \Scal_i^\Abf\cup \Fcal_i^\Abf \text{ and $i$ dominant,} 
\ee 
which implies \eqref{contraintes-i2}. 
\end{proof}

\begin{proof}[Proof of Corollary~\ref{corol}]  
From \eqref{matrixp} we see that the eigenvectors of the averaged matrix associated with the $p$-system
$$
\rb_j(v,v') = \begin{pmatrix} 
\pm \cb(v,v')
\\
1
\end{pmatrix},
$$ 
only depend upon the eigenvalues $\pm \cb(v,v')$. This shows that \eqref{eigenv2} holds for the $p$-system.

Next, we turn attention to the Euler equations of isentropic fluid dynamics, 
\be
\label{Euler}
\del_t \rho + \del_x (\rho u) = 0, \qquad
\del_t(\rho u) + \del_x (\rho u^2 + p(\rho)) = 0, 
\ee 
where $\rho \geq 0$ and $u$ represents the specific density and velocity of the fluid, respectively. The pressure
$p=p(\rho)$ is assumed to be strictly increasing for $\rho>0$, so that the equations are strictly hyperbolic, 
away from the vacuum ($\rho=0$) at least. 
We use here the density $\rho$ and the momentum $q : = \rho u$, which are the conservative variables of \eqref{Euler}.  
To define the averaged matrix we consider two solutions $(\rho,q)$, $(\rho',q')$ of \eqref{Euler} 
and write  
$$
\del_t (\rho-\rho') + \del_x (q-q' ) = 0, 
\qquad
\del_t(q-q') + \del_x \Big({q^2\over \rho} - {{q'}^2 \over \rho'} + p(\rho) - p(\rho')\Big) = 0,
$$
which we now transform into the form \eqref{1.5}. 

Let us introduce a function $\eb=\eb(\rho,q,\rho',q')$ satisfying 
$$
{q^2\over \rho} - {{q'}^2 \over \rho'} = - \eb^2 (\rho - \rho') + 2 \eb (q-q'), 
$$
that is, 
\be
\label{ebar}
\aligned 
\eb(\rho,q,\rho',q') 
& := {\rho u  - \rho' u' - \sqrt{\rho \rho'} (u - u') \over \rho - \rho'}
\\
& = {u + u' \over 2} + {1 \over 2} {\sqrt{\rho} - \sqrt{\rho'} \over \sqrt{\rho} + \sqrt{\rho'}} (u - u'). 
\endaligned 
\ee
Introduce also the function $\cb:\R^2 \to \R$ by 
\be
\label{cbar}
\cb^2(\rho,\rho') := {p(\rho) - p(\rho') \over \rho - \rho'}. 
\ee

Then, the maps $\varphi := \rho - \rho'$ and $\psi := q - q'$ satisfy the linear hyperbolic system
$$
\del_t \varphi + \del_x \psi = 0, 
\qquad
\del_t \psi + \del_x (- \eb^2 \varphi + 2\eb \psi + \cb^2(\rho,\rho') \, \psi) = 0, 
$$ 
which corresponds to the following choice of averaged matrix:  
\be
\label{matrixE}
\Ab(\rho,q; \rho';q') := \begin{pmatrix} 0 & 1
\\
-\eb^2 + \cb^2  &  2 \eb
\end{pmatrix}. 
\ee
Clearly, the eigenvalues of the matrix $\Ab(\rho,q; \rho';q')$ are $\lam_j = \eb \pm \cb$, 
while the eigenvectors are $r_j = (1, \eb \pm \cb)^T$. This establishes the desired property \eqref{eigenv2}
for the Euler equations. 

Consequently, Theorem~\ref{theorem-main} applies to both the Lagrangian and the Eulerian formulations of the fluid
dynamics equations, and this establishes \eqref{5.4} and \eqref{5.5} for the fluid dynamics systems. 
\end{proof}


\section{Concluding remark} 

The construction of the weight proposed in Section~\ref{WE-0} can be generalized, 
by noticing that the jump of the weight need not be related to the averaged speed $\ab$.

\begin{proposition}
\label{7.1} 
Let $\pi : \R^2 \to \R$ be a given smooth function. 
Under the assumptions of Theorem~\ref{51} there exists a weight $\wbf$ satisfying, instead
of \eqref{contr-scalar}, the more general condition
\be
\label{contr-scalar44}
\wbf^+ - \wbf^-   = 
\begin{cases}  
- K \,  (\pi^+ - \pi^-),  & (\psi ^+ - \psi^-) \, \psi^+ >0 \text { and } \psi^- \psi^+ >0.   
\\ 
\hskip.3cm  K \, (\pi^+ - \pi^-),  & (\psi ^+ - \psi^-) \, \psi^+ < 0 \text { and } \psi^- \psi^+ >0, 
\end{cases}
\ee
where $\psi := u_h - u_h'$ and $\pi^\pm =\pi(u_\pm^h, {u_\pm^h}')$.
\end{proposition}

The proof is immediate from the one of Theorem~\ref{51}. Of course, we recover the previous setting
by choosing $\pi = \ab$. 

As far as the application to hyperbolic conservation laws is concerned and 
in order to use this weight in the weighted norm however, we need \eqref{updown} to hold
at undercompressive waves, that is in the region where $\psi$ keeps a constant sign. 
This essentially forces us to pick up $\pi= \ab$ in Proposition~\ref{7.1}.


\section*{Acknowledgements}  

The author was partially supported by the University of Cambridge (UK) 
and the Centre National de la Recherche Scientifique (CNRS), and 
also gratefully acknowledges the support and hospitality of the Isaac Newton Institute 
for Mathematical Sciences, during the Semester Program 
{\sl ``Nonlinear Hyperbolic Waves in Phase Dynamics and Astrophysics'',}
organized by C.M. Dafermos, P.G. LeFloch, and E. Toro.  This work was completed at the Mittag-Leffler Institute 
(Stockholm) during the Semester Program {\sl ``Nonlinear Wave Motion'',}  
organized by A. Constantin, C.M. Dafermos, H. Holden, K.H. Karlsen, and W. Strauss.


\end{document}